\documentclass[12pt]{article}
\usepackage{latexsym}
\usepackage{amsfonts}
\usepackage{amssymb}
\date{}
\newtheorem{Th}{\quad Theorem}[section]
\newtheorem{Lm}{\quad Lemma}[section]
\newtheorem{Prop}{\quad Proposition}[section]

\title{A complete rewrite system and normal forms for $(S)_{\mathrm{reg}}$ }

\author{Jean-Camille Birget \\  
       {\footnotesize Dept.\ of Computer Science} \\
        {\footnotesize Rutgers University -- Camden} \\
        {\footnotesize Camden, NJ 08102, USA} \\
        {\footnotesize birget@camden.rutgers.edu }
\and
   Stuart W. Margolis\thanks{Both authors' research was supported in part by 
NSF grant DMS-9203981} \\ 
        {\footnotesize Dept.\ of Mathematics and Computer Science} \\
        {\footnotesize Bar Ilan University} \\
        {\footnotesize Ramat Gan 52900, Israel} \\
        {\footnotesize margolis@macs.biu.ac.il }
        }

\begin{document}

\maketitle

\noindent {\bf Abstract.} \  The $(.)_{\mathrm{reg}}$ construction was 
introduced in order to make an arbitrary semigroup $S$ divide a regular 
semigroup $(S)_{\mathrm{reg}}$ which shares some important properties with 
$S$ (e.g., finiteness, subgroups, torsion bounds, $J$-order structure). 
We show that $(S)_{\mathrm{reg}}$ can be described by a rather simple 
{\it complete string rewrite system}, 
as a consequence of which we obtain a new proof of the normal 
form theorem for $(S)_{\mathrm{reg}}$. 
The new proof of the normal form theorem is conceptually simpler than the
previous proofs.

\section{Introduction}
 
Regular semigroups have always played a special role in the structure theory 
of semigroups. Since, however, semigroups are in general not regular, it is 
interesting to connect arbitrary semigroups to regular ones. An obvious 
connection of this sort is the embedding of any semigroup $S$ into a full 
transformation semigroup (which is always a regular semigroup). 
A much tighter connection was proved in \cite{Bi1}, \cite{Bi2}: Any 
semigroup $S$ divides a regular semigroup $(\hat{S})_{\mathrm{reg}}$; if $S$ 
is finite, then $(\hat{S})_{\mathrm{reg}}$ is finite; every subgroup of 
$(\hat{S})_{\mathrm{reg}}$ divides a subgroup of $S$; 
$(\hat{S})_{\mathrm{reg}}$ has the same regular ${\cal J}$-order as $S$, and shares 
many other properties with $S$. 

In more detail, the division of $S$ into $(\hat{S})_{\mathrm{reg}}$ is done 
in two steps: First $S$ is expanded to the left-right-iterated Rhodes
expansion $\hat{S}$; this yields an {\it unambiguous} semigroup, i.e., a 
semigroup whose  ${\cal L}$-order and ${\cal R}$-order are forests  \cite{Bi1}.
Then $\hat{S}$ is embedded into the regular semigroup $(\hat{S})_{\mathrm{reg}}$
by applying the $(.)_{\mathrm{reg}}$ construction \cite{Bi2}. When $S$ is any
unambiguous semigroup then $S$ is a subsemigroup of $(S)_{\mathrm{reg}}$;
when $S$ is not unambigous then $S$ is not a subsemigroup of 
$(S)_{\mathrm{reg}}$; in that case, the subsemigroup of 
$(S)_{\mathrm{reg}}$ generated by $S$ is the Rees quotient of $S$ over 
the ideal of {\it ambiguous} elements of $S$ (by definition, an element $s \in S$
is ambiguous iff the ${\cal L}$-order and ${\cal R}$-order above $s$ are not both
forests \cite{Bi2}). 

As a consequence of this, 
{\it  every aperiodic (finite) semigroup divides a regular 
aperiodic (finite) semigroup.}
By definition, a semigroup is aperiodic iff it 
satisfies the identity $x^n = x^{n+1}$ for some positive integer $n$.
More generally, 
{\it  an infinite torsion semigroup (or a bounded torsion semigroup,
satisfying $x^t = x^{t+c}$) divides a regular torsion semigroup (respectively, 
a bounded torsion semigroup satisfying $x^{t+1} = x^{(t+1)+c}$).
Also, a semigroup whose subgroups belong to some variety (or pseudo-variety, 
or quasi-variety) ${\cal G}$ divides a regular semigroup whose 
subgroups belong to the same variety (resp., pseudo-variety, or quasi-variety)
 }  ${\cal G}$. 
So far, the above method is the only known proof of these results.

Another application of $(.)_{\mathrm{reg}}$ is to find an improved version of the 
Rhodes-Allen Synthesis theorem, a generalization of both the Rees theorem and
the Krohn-Rhodes theorem (see \cite{Bi4} and \cite{GrBook} for background). 

The $(.)_{\mathrm{reg}}$ construction itself has connections with two-way finite 
automata \cite{Pecuchet}.

As we will see below, the $(.)_{\mathrm{reg}}$ construction is rather easy to 
describe, but it is not easy to prove {\it the normal form theorem} for the 
elements of $(S)_{\mathrm{reg}}$. The normal form is important,
because it is used to prove the main properties of $(S)_{\mathrm{reg}}$.
However, the fact that $S$ is a 
subsemigroup of $(S)_{\mathrm{reg}}$ when $S$ is unambiguous, has a relatively
simple direct proof -- see \cite{Bi2}, pp.\ 73-75.
All the known proofs of the normal form theorem are tedious. 
The original proof of the normal form theorem 
in \cite{Bi2} uses Van der Waerden's method (letting $S$ act faithfully on a
set of normal forms). More recently, Grillet \cite{GrPaper} introduced
another method, based on congruences on non-associative structures.
The present paper contains a third proof, based on string rewriting. Besides 
providing yet another proof, we show that $(S)_{\mathrm{reg}}$ can be defined by a
rather simple {\em complete string rewrite system}; this makes the normal forms of the 
elements of $(S)_{\mathrm{reg}}$ obvious. Unfortunately, the catch is that the 
confluence of this rewrite system requires a rather tedious proof,
though, conceptually, this proof is rather easy and looks almost like a verification 
by a machine. 

We will assume from now on that $S$ is unambiguous.

\bigskip

\noindent {\bf Notation and definitions:} 

By $>_{\cal L}, \ \leq_{\cal L}, \ \equiv_{\cal L}$ we denote
Green's well known ${\cal L}$-relations, and similarly for the 
${\cal R}$-relations. 
We also use the ${\cal D}$-equivalence \ $\equiv_{\cal D}$. 
See e.g. \cite{GrBook} for background.

We will also need the ${\cal L}$-{\em incomparability} relation
\ $\frac{<}{>}\!\!\!|\,_{_{\cal L}}$ \ defined as follows: \    
$s \ \frac{<}{>}\!\!\!|\,_{_{\cal L}} \ t$ \ iff neither $s \leq_{\cal L} t$
nor $s \geq_{\cal L} t$. We also define ${\cal L}$-{\em comparability}: \ 
$s \ \frac{<}{>}\!\!\,_{_{\cal L}} \ t$ \ iff either $s \leq_{\cal L} t$
or $s \geq_{\cal L} t$. \ A similar notation is used for ${\cal R}$. 

Following \cite{Bi1}, \cite{Bi2}, we call a semigroup $S$ {\em unambiguous} 
iff \ for all $s, t, u \in S - \{ 0 \}$ : \ 
$s >_{\cal L} u <_{\cal L} t$ \ implies 
\ $s \ \frac{<}{>}\!\!\,_{_{\cal L}} \ t$ \, and 
\ $s >_{\cal R} u <_{\cal R} t$ \ implies 
\ $s \ \frac{<}{>}\!\!\,_{_{\cal R}} \ t$.  \  This means that the ${\cal L}$-order 
on the ${\cal L}$-classes of $S - \{ 0 \}$ is a forest, and similarly for
${\cal R}$.    
(Here, 0 is the zero of $S$ if $S$ has a zero; otherwise, $S - \{ 0 \} = S$.)

In order to avoid confusion between products of elements in a semigroup $S$ 
and strings of elements of $S$, we denote a {\em string} of length $n$ as 
an $n$-tuple of the form $(s_1, s_2, \ldots , s_n)$. The product of these 
elements in $S$ is denoted by \ $s_1 s_2 \ldots s_n$ \ or \ 
$s_1 \cdot s_2 \cdot \ldots \cdot s_n$ ($\in S$).

When $S$ does not have an identity element, $S^{1}$ denotes the monoid obtained
by adding a new identity element to $S$; if $S$ is already a monoid, 
$S^{1}$ is just $S$.

We refer to \cite{Ja} for background on rewrite systems.

\section{The rewrite system}

\noindent {\bf Presentation of  $(S)_{\mathrm{reg}}$ by generators and relations:}

\medskip

Let $S$ be a semigroup (possibly infinite). Let 0 be the zero of $S$, 
if $S$ has a zero; otherwise, let 0 be a new symbol not in $S$.
Let $\overline{S - \{0\}} = \{ \overline{s} : s \in S - \{0\} \}$ be a set 
that is disjoint from $S \cup  \{0\}$, where the map \ 
$x \in (S - \{0\}) \cup \overline{S - \{0\}} \longmapsto \overline{x} \in 
(S - \{0\}) \cup \overline{S - \{0\}}$ \ 
is a bijection such that \ $\overline{\overline{x}} = x$. 
We also let $\overline{0} = 0$; the symbol $\overline{0}$ will never be used 
and will always automatically be replaced by $0$. 

Then, following \cite{Bi2}, $(S)_{\mathrm{reg}}$ is defined by the following 
presentation:

\medskip

\noindent {\it Generators}:  

\smallskip

$S \ \cup \ \overline{S - \{0\}} \ \cup \ \{0\}$.

\medskip

\noindent {\it Relations}:

\smallskip

$(s, t) = (s t)$ \ \ \ \  for all $s, t \in S$  

\smallskip

$(\overline{s}, \overline{t}) = (\overline{t s})$ \ \ \ \ for all 
$s, t \in S - \{0\}$ 

\smallskip

$(0) = (0, 0) = (0, s) = (s, 0) = (0, \overline{s}) = (\overline{s}, 0)$ 
\ \ \ \ for all $s \in S$

\smallskip

$(s, \overline{t}) = (0)$ \ \ \ if \  $s \ \frac{<}{>}\!\!\!|\,_{_{\cal L}}
 \ t$, \ \ \ $s, t \in S - \{0\}$

\smallskip

$(\overline{s}, t) = (0)$ \ \ \ if \ $s \ \frac{<}{>}\!\!\!|\,_{_{\cal R}} \ 
t$, \ \ \ $s, t \in S - \{0\}$

\smallskip

$(s, \overline{s}, s) = (s)$ \ \ \ \  for all $s \in S - \{0\}$

\smallskip

$(\overline{s}, s, \overline{s}) = (\overline{s})$ \ \ \ \ 
for all $s \in S - \{0\}$

\medskip

It is proved in \cite{Bi2} (see also \cite{GrPaper} and \cite{GrBook}) that 
$S$ is a subsemigroup of $(S)_{\mathrm{reg}}$ if $S$ is unambiguous, and that
$(S)_{\mathrm{reg}}$ is a regular semigroup with involution (i.e., for all
$x, y \in (S)_{\mathrm{reg}}: \ \overline{\overline{x}} = x, \ 
\overline{xy} = \overline{y} \, \overline{x}, \ x \, \overline{x} \, x = x)$.

\bigskip

\noindent {\bf Rewrite rules for $(S)_{\mathrm{reg}}$:}

\medskip

We now introduce a string rewrite system for $(S)_{\mathrm{reg}}$. This rewrite
system is finite iff $S$ is finite.
The reduced words of this rewrite system are the normal forms of $(S)_{\mathrm{reg}}$.
In the next sections we will prove that this rewrite system is complete, when
$S$ is unambiguous.

\medskip

\noindent {\bf 1.} {\em Length-reducing rules:}

\smallskip 

The last two of the following set of rules make use of a partial function \ 
$B : S \times S \times S \to S$, \ that will be defined after the statement 
of all the rules. 

\bigskip

\noindent (1.1) \ \ \ \ \ \ $(s, t) \to (s t)$, 
 \ $(\overline{s}, \overline{t}) \to (\overline{t s})$
\ \ \ \ for all $s, t \in S$

\smallskip

\noindent (1.2) \ \ \ \ \ \ $(0, 0) \to (0), \ (0, s) \to (0), 
\ (0, \overline{s}) \to (0), \ (\overline{s}, 0) \to (0)$ 
\ \ \ \ for all $s \in S$

\medskip

\noindent (1.3) \ \ \ \ \ \ $(s, \overline{t}) \to (0)$ \ \ \ \ if \ 
$s \ \frac{<}{>}\!\!\!|\,_{_{\cal L}} \ t$, \ \ \ $s, t \in S -\{0\}$

\medskip

\noindent (1.4) \ \ \ \ \ \ $(\overline{s}, t) \to (0)$ \ \ \ \ if \ 
$s \ \frac{<}{>}\!\!\!|\,_{_{\cal R}} t$, \ \ \ $s, t \in S -\{0\}$

\medskip

\noindent (1.5) \ \ \ \ \ \ $(u, \overline{v}, w) \to (B(u,v,w))$ \ \ \ if \ 
$u \leq_{_{\cal L}} v \geq_{_{\cal R}} w$, \ \ \ $u, v, w \in S -\{0\}$

\medskip

\noindent (1.6) \ \ \ \  \ \ $(\overline{u}, v, \overline{w}) \to 
(\overline{B(w,v,u)})$ \ \ \ \ 
if \ $u \leq_{_{\cal R}} v \geq_{_{\cal L}} w$, \ \ \ $u, v, w \in S -\{0\}$

\bigskip

\noindent {\bf 2.} {\em Length-preserving rules:}

\smallskip

For these rules we choose one {\it representative} element in every 
${\cal R}$-class and in every ${\cal L}$-class. We make these choices 
so that ${\cal D}$-related representatives of ${\cal R}$-classes 
are ${\cal L}$-related, and 
${\cal D}$-related representatives of ${\cal L}$-classes
are ${\cal R}$-related. Moreover, if two representatives (one representing
an ${\cal L}$-class and one representing an ${\cal R}$-class)
are in the same ${\cal H}$-class they are chosen to be equal.
Such a choice can always be made.

Note that this condition on the choice of representatives was not used, and 
not required, in \cite{Bi2} and \cite{GrPaper}. A similar choice however is 
made in the Rees-Sushkevitch coordinatization, see e.g. \cite{GrBook}. 
Notation: For any $s \in S$ the chosen representative of the 
${\cal R}$-class (or ${\cal L}$-class) of $s$ is $r_s$ 
(respectively $\ell_s$).

The length-preserving rules make use of two partial functions, 
$B_{_{\cal L}}$ and $B_{_{\cal R}} :
S \times S \to S$, that will be defined after the rules.

Note the unsymmetry between rules (2.1)-(2.2) and (2.3)-(2.4), which
is needed for obtaining unique normal forms; see \cite{Bi2}, \cite{GrPaper}
for more discussion on the normal forms.  

\bigskip

\noindent (2.1) \ \ \ \ \ \ $(s, \overline{t}) \to 
(r_s, \overline{B_{_{\cal R}}(s,t)})$
\ \ \ \ if \ $s >_{_{\cal L}} t$ \ and \ $s \neq r_s$ 

\medskip

\noindent (2.2) \ \ \ \ \ \ $(\overline{s}, t) \to (\overline{\ell_s}, 
B_{_{\cal L}}(t,s))$ \ \ \ \ if \ $s >_{_{\cal R}} t$ \ and \ $s \neq \ell_s$

\medskip

\noindent (2.3) \ \ \ \ \ \ $(\overline{t}, s) 
\to (\overline{B_{_{\cal L}}(t,s)}, \ell_s)$ \ \ \ \ 
if \ $t \leq_{_{\cal R}} s$ \ and \ $s \neq \ell_s$

\medskip

\noindent (2.4) \ \ \ \ \ \ $(t, \overline{s})  \to (B_{_{\cal R}}(s,t), 
\overline{r_s})$ \ \  
\ \ if \ $t \leq_{_{\cal L}} s$ \ and \ $s \neq r_s$

\bigskip

\noindent {\bf Definition of $B$.} \  {\it 
If \ $u \leq_{_{\cal L}} v \geq_{_{\cal R}} w$, where
$u, v, w \in S - \{0\}$, then \ $B(u, v, w) = u z$, where $z \in S^{1}$ 
is such that $w = v z$. }

\smallskip

This  operation was used in  \cite{Bi2}, but was first explicitly defined in
\cite{GrPaper}.  It is easy to see that if $u \leq_{_{\cal L}} v 
\geq_{_{\cal R}} w$
then $B(u, v, w)$ exists and is unique (i.e., it depends only on $u, v, w$ and 
not on $x$; see Lemma \ref{Lm3} below). The main motivation for $B$ is that
in $S_{\rm {reg}}$, \ $u \overline{v} w = B(u, v, w)$ \ if 
$u \leq_{_{\cal L}} v \geq_{_{\cal R}} w$, as we will prove in Proposition 2.1 
below. 

\bigskip

\noindent {\bf Definition of $B_{_{\cal R}}$ and $B_{_{\cal  L}}$.} \  
{\it If \ $u \geq_{_{\cal L}} v$, \ where $u, v \in S - \{0\}$, then \  
$B_{_{\cal R}}(u, v) = x r_u$, where $x \in S^{1}$ is such that $v = x u$. 
\ If \ $v \leq_{_{\cal R}} u$,  where $u, v \in S - \{0\}$, then \  
$B_{_{\cal  L}}(v, u) = \ell_u y$, where $y \in S^{1}$ is such that
$v = u y$. }

\smallskip

This  operation was implicit in \cite{Bi2}. Again, it is easy to see that if 
$u \geq_{_{\cal L}} v$ (or $v \leq_{_{\cal R}} u$) then $B_{_{\cal R}}(u, v)$ 
(resp. $B_{_{\cal  L}}(v, u)$) exists and is unique (i.e., it depends only
on $u$ and $v$). The main motivation for $B_{_{\cal R}}$ is that
in $S_{\rm {reg}}$, \ $v \overline{u} = B_{_{\cal R}}(u,v) \, \overline{r_u}$
 \  if $u \geq_{_{\cal L}} v$, as we will prove in Proposition 2.1 below. 
The motivation for $B_{_{\cal  L}}$ is similar.

In the next section we will see another, pictorial motivation for
$B$, $B_{_{\cal R}}$, and $B_{_{\cal  L}}$.

\smallskip

Before proving the next proposition we need to recall a key property of
$(S)_{\mathrm{reg}}$.

\begin{Lm} {\rm (Fact 2.5 in \cite{Bi2}).}
For all $s, r, \ell \in S - \{0\}$: 
If $s \equiv_{_{\cal R}} r$ in $S$ then $s \, \overline{s} = r \, \overline{r}$ 
in $(S)_{\mathrm{reg}}$. \  
If $s \equiv_{_{\cal L}} \ell$ in $S$ then 
$\overline{s} \, s = \overline{\ell} \, \ell$ in $(S)_{\mathrm{reg}}$.
\end{Lm}
{\it Proof.} Let $a, b \in S^{1}$ be such that $s = ra, \ r = sb$; so,
$rab = r$. Then, by using the relations of the presentation of 
$(S)_{\mathrm{reg}}$ we have:

$s \, \overline{s} \ = \ ra \ \overline{ra} \ = \ 
r \, \overline{r} \, r \ a \ \overline{ra} \ = \   
r \ \overline{rab} \ ra \ \overline{ra} \ = \   
r \, \overline{b} \ \overline{ra} \ ra \ \overline{ra} \ =  \   
r \, \overline{b} \ \overline{ra} \ = \ r \ \overline{rab} \ = \   
r \, \overline{r}$.  
\ \ \  $\Box$

\begin{Prop} 
The rewrite system defines $(S)_{\mathrm{reg}}$.
\end{Prop}
{\it Proof.}  The rewrite rules (when made symmetric) imply the relations of 
the presentation; to obtain the last two relations of the presentation, let 
$u = v = w$ in rules (1.5) and (1.6).

Conversely it is straightforward to show that in $(S)_{\mathrm{reg}}$ the 
relations corresponding to the rules (1.5), (1.6), (2.1)--(2.4) hold (see 
also \cite{Bi2}). 

\smallskip

Let us derive rule (1.5). Since $u \leq_{_{\cal L}} v \geq_{_{\cal R}} w$, 
let $x, y \in S^{1}$ be such that $u = xv, \ w = vy$. Then \  
$u \, \overline{v} \, w \ = \ x \, v \, \overline{v} \, v \, y  \ = \ xvy$, 
using $v \overline{v} v = v$ in $(S)_{\mathrm{reg}}$. Moreovre,
$xvy  = uy = B(u,v,w)$ by the definition of $B$. Thus, \  
$u \, \overline{v} \, w \ = \ B(u,v,w)$ \ in $(S)_{\mathrm{reg}}$.

\smallskip

Let us derive rule (2.4). Since $t \leq_{_{\cal L}} s$, let $x \in S^{1}$ be 
such that $t = xs$. Then $t \overline{s} \ = \ x \, s \, \overline{s} \ = \ 
x \, r_s \, \overline{r_s}$; the last equality follows from the last Lemma.
And $x r_s = B_{\cal R}(s,t)$, by the definition of $B_{\cal R}$.
Thus, \ $t \overline{s} =  B_{\cal R}(s,t) \, \overline{r_s}$
 \ in $(S)_{\mathrm{reg}}$.  

\smallskip
 
The other rules can be derived in a very similar way. \ \ \  $\Box$

\medskip

\noindent
One of the main results of \cite{Bi2} is the following:

\smallskip

\noindent
{\bf Normal Form theorem for} $(S)_{\mathrm{reg}}$: \ 
If $S$ is {\em unambiguous} then $S$ is a subsemigroup of $(S)_{\mathrm{reg}}$, and
(for any fixed choice of representatives of the ${\cal L}$- and 
${\cal R}$-classes)
every element of $(S)_{\mathrm{reg}}$ can be written in a unique way in the 
{\bf normal form} 

\medskip

(0) \ \  or

\medskip

$([r_1,] \ \overline{\ell_2}, \ \ldots \ , r_{n-1}, \, \overline{\ell_n}, \ s, \ 
\overline{r'_m}, \, \ell'_{m-1}, \  \ldots \ , \overline{r'_2} \ [, \ell'_1])$

\medskip

\noindent where \ \ \ 
$[r_1 >_{\cal L}] \  \ell_2 >_{\cal R} \ \ldots \ >_{\cal R} r_{n-1} 
>_{\cal L} \ell_n 
>_{\cal R} s \leq_{\cal L} r'_m <_{\cal R} \ell'_{m-1} <_{\cal L} \ \ldots \  
<_{\cal L} r'_2 \ [ <_{\cal R} \ \ell'_1]$, \ 

\medskip

\noindent or in the form 

\medskip

$([r_1,] \ \overline{\ell_2}, \ \ldots \ , r_{n-2}, \, \overline{\ell_{n-1}}, \,  r_n, 
\ \overline{s}, \ \ell'_m, \, \overline{r'_{m-1}}, \, \ell'_{m-2}, \  \dots \ ,  
\overline{r'_2} \ [, \ell'_1])$

\medskip

\noindent where \ \ \ 
$[r_1 >_{\cal L} ] \  \ell_2 >_{\cal R} \ \ldots \ >_{\cal R} r_{n-2} >_{\cal L}
\ell_{n-1} >_{\cal R} r_n >_{\cal L} s \leq_{\cal R} \ell'_m <_{\cal L} 
r'_{m-1} <_{\cal R} \ell'_{m-2} <_{\cal L} \ \ldots \ 
<_{\cal L} r'_2 \ [ <_{\cal R} \ \ell'_1]$. 

\medskip

\noindent Here, every $r_i, r'_j, \ell_i, \ell'_j$ is a representative of an 
${\cal R}$- or ${\cal L}$-class, and $s$ is any element of $S - \{0\}$.
Elements in square brackets may be absent.

The normal form representation is the key to many structure properties of 
$(S)_{\mathrm{reg}}$, e.g., the fact that $S$ and $(S)_{\mathrm{reg}}$ have the same 
${\cal J}$-class structure. 
The main result of this paper is:

\begin{Th}
The above rewrite system for $(S)_{\mathrm{reg}}$ is complete (i.e., confluent and
terminating). The normal forms of the rewrite systems are as given above.
\end{Th}

The remainder of this paper consists of the proof of this theorem.
In Section 3 we give some basic properties of  $B$, $B_{_{\cal  L}}$, 
and $B_{_{\cal R}}$, then in Section 4 we prove termination of the 
rewrite system, and finally in Section 5 we prove local confluence.

\section{Properties of the functions $B$, $B_{_{\cal  L}}$, and $B_{_{\cal R}}$}

In this section we collect all the basic properties of 
$B$, $B_{_{\cal  L}}$, and $B_{_{\cal R}}$ that we will need in order to 
prove that
the rewrite system for $(S)_{\mathrm{reg}}$ is terminating and locally confluent. 
The reader may skip this section, and come back to it while reading the proofs
of termination and local confluence. 

Below, when we write an expression like $B_{_{\cal R}}(x, y)$, 
$B_{_{\cal L}}(x, y)$, or $B(x, y, z)$, we  always implicitly assume that these 
expressions are defined (i.e., we assume that $x \geq_{_{\cal L}} y$ when 
we use $B_{_{\cal R}}(x, y)$, etc.).

\bigskip

In all the proofs in this section it will be useful for the reader to 
represent $B$, $B_{_{\cal R}}$, and $B_{_{\cal L}}$, by the following diagrams, 
which are justified by the next few lemmas. 

\bigskip

\noindent {\bf Diagram of} $B(u,v,w)$:

\medskip

If \ $u \leq_{_{\cal L}} v \geq_{_{\cal R}} w$, let $y, z \in S^{1}$ be any 
elements such that $u = yv$, $w = vz$. Then we have the commutative diagram:

\bigskip

\unitlength=0.90mm \special{em:linewidth 0.4pt}
\linethickness{0.4pt}
\begin{picture}(139,50)

\put(34,43){\makebox(0,0)[cc]{$v$}}             
\put(18,26){\makebox(0,0)[cc]{$u$}}             
\put(50,26){\makebox(0,0)[cc]{$w$}}           
\put(34,8){\makebox(0,0)[cc]{$B(u,v,w)$}}     

\put(30,40){\vector(-1,-1){10}}  
\put(38,40){\vector(1,-1){10}}  
\put(21,22){\vector(1,-1){10}}  
\put(47,22){\vector(-1,-1){10}}  

\put(23,37){\makebox(0,0)[cc]{$y \cdot$}}
\put(45,37){\makebox(0,0)[cc]{$\cdot z$}}

\put(22,17){\makebox(0,0)[cc]{$\cdot z$}}
\put(47,17){\makebox(0,0)[cc]{$y \cdot$}}

\end{picture}

\bigskip

\bigskip

\noindent {\bf Diagram of} $B_{_{\cal R}}(u,v)$:

\medskip

If \ $s \geq_{_{\cal L}} t$, let $x \in S^{1}$ be any
elements such that $t = xs$. Also, let $a,a' \in S^{1}$ be such that 
$r_s a = s$ ans $r_s = s a'$. Then we have the commutative diagram:

\bigskip

\unitlength=0.90mm \special{em:linewidth 0.4pt}
\linethickness{0.4pt}
\begin{picture}(139,50)

\put(18,38){\makebox(0,0)[cc]{$s$}}             
\put(18,8){\makebox(0,0)[cc]{$t$}}              
\put(50,38){\makebox(0,0)[cc]{$r_s$}}           
\put(54,8){\makebox(0,0)[cc]{$B_{\cal R}(s,t)$}}     

\put(24,39){\vector(1,0){20}}  
\put(44,37){\vector(-1,0){20}}  
\put(24,9){\vector(1,0){20}}  
\put(44,7){\vector(-1,0){20}}  

\put(34,42){\makebox(0,0)[cc]{$\cdot a'$}}             
\put(34,35){\makebox(0,0)[cc]{$\cdot a$}}             
\put(34,12){\makebox(0,0)[cc]{$\cdot a'$}}             
\put(34,5){\makebox(0,0)[cc]{$\cdot a$}}             

\put(18,34){\vector(0,-1){22}} 
\put(50,34){\vector(0,-1){22}} 
\put(15,24){\makebox(0,0)[cc]{$x \cdot$}}
\put(53,24){\makebox(0,0)[cc]{$x \cdot$}}
\end{picture}

\bigskip

The diagram for $B_{_{\cal L}}$ is similar to the diagram for $B_{_{\cal R}}$.

\bigskip

\begin{Lm} \label{Lm1} \ 
{\bf (a)} \ If \  $u = r_u \alpha$ \  then \ 
$B_{_{\cal R}}(u, v) \cdot \alpha = v$.
\ Similarly, if \ $v = \beta \ell_u$  then $\beta \cdot B_{_{\cal  L}}(v, u) 
= u$. \\ 
{\bf (b)} \ If \  $r_u = u \alpha'$ \ then \ $B_{_{\cal R}}(u, v) = v \alpha'$. \ 
Similarly, if \ $\ell_u = \beta' u$ then \ $\beta \cdot B_{_{\cal L}}(v, u) = 
\beta' u$.
\end{Lm}
The proof is trivial.

\begin{Lm} \label{Lm2}
$B_{_{\cal R}}(r_u, v) = v$, and \ $B_{_{\cal L}}(v, \ell_u) = v$.
\end{Lm}
The proof is trivial.

\begin{Lm} \label{Lm3}
$B_{_{\cal R}}(u, v) \equiv_{_{\cal R}} v$, and $B_{_{\cal L}}(v, u) 
\equiv_{_{\cal L}} v$. 
\end{Lm}
{\it Proof.} \  If we multiply $r_u  \equiv_{_{\cal R}} u$ on the left by 
$x$ we obtain $B_{_{\cal R}}(u, v) = x  r_u \equiv_{_{\cal R}} x u = v$.
For ${\cal L}$ the proof is similar. 
\ \ \ $\Box$

\begin{Lm} \label{Lm4}
If $s <_{_{\cal L}} t$ then \ $B_{_{\cal R}}(t, s) <_{_{\cal L}} r_t$ \ 
(and the same holds with
$<_{_{\cal L}}$ replaced by $\equiv_{_{\cal L}}$ or $\leq_{_{\cal L}}$). \  
If $t >_{_{\cal R}} s$ then \ $\ell_t >_{_{\cal R}} B_{_{\cal L}}(s, t)$ \  
(and the same holds
with $<_{_{\cal R}}$ replaced by $\equiv_{_{\cal R}}$ or $\leq_{_{\cal R}}$).
\end{Lm}
{\it Proof.} \ We prove the first statement, the other ones having very similar
proofs. Let $a$ be such that $ta = r_t$.

Since $B_{_{\cal R}}(t, s) = xr_t$ for some $x$ such that $xt = s$,
we have $B_{_{\cal R}}(t, s) = xr_t \leq_{_{\cal L}} r_t$. Actually we have 
$B_{_{\cal R}}(t, s) <_{_{\cal L}} r_t$. Indeed, if we had 
$xr_t \equiv_{_{\cal L}} r_t$, then multiplying on the right by $a$ yields 
$s = xr_ta \equiv_{_{\cal L}} r_ta = t$,
i.e., $s \equiv_{_{\cal L}} t$, which contradicts the assumption.
\ \ $\Box$

\begin{Lm} \label{Lm5}
If \ $u \leq_{\cal L} v \geq_{\cal R} w$ \ then \ 
$B(u, v, w) = yw = ux = yvx$, \ 
where $x$ is such that $w = vx$, and $y$ is such that $u = yv$. The value of
$B(u, v, w)$ does not depend on the $x$ or $y$ chosen.
\end{Lm}
{\it Proof.} \ By definition, $B(u, v, w) = ux$ where $x$ is such that 
$w = vx$.  Hence $B(u, v, w) = ux = yvx = yw$. 

To see that $B(u, v, w)$ does not depend on the choice of $x$ (provided that
$w = vx$), let $w = vx_1 = vx_2$. Then $B(u, v, w) = yvx_1 = yvx_2$.  
Similarly, one sees that the choice of $y$ does not matter (provided that 
$u = yv$).  \ \ $\Box$

\begin{Lm} \label{Lm6}
If \ $u \leq_{\cal L} v \geq_{\cal R} w$ \ and $t \in S - \{0\}$ then \  
$B(tu, v, w) = t \cdot B(u, v, w)$ \ and \  $B(u, v, wt) = B(u, v, w) \cdot t$.
\end{Lm}
{\it Proof.} \ Since $B(u, v, w) = ux$ where $x$ is such that $w = vx$, we
obtain $t \cdot B(u, v, w) = tux$ with $w = vx$. Hence by the definition of
$B(tu, v, w)$ we have $B(tu, v, w) = t \cdot B(u, v, w)$. 

The proof for $B(u, v, wt)$ is similar, by using Lemma \ref{Lm5}.
\ \ $\Box$ 

\begin{Lm} \label{Lm7}  \ \  \\ 
(1) \ If \ $u \geq_{\cal L} su \geq_{\cal L} v$ \ then \ 
$sr_u \geq_{\cal L} B_{\cal R}(u, v)$ \ and \   
$B_{\cal R}(su, v) = B_{\cal R}(sr_u, B_{\cal R}(u, v))$. \\
(2) \ If \ $su \leq_{\cal L} v \leq_{\cal L} u$ \ then \  
$sr_u \leq_{\cal L} B_{\cal R}(u, v)$ \ and \  
$B_{\cal R}(v, su) = B_{\cal R}(B_{\cal R}(u, v), sr_u)$. \\
If \ $su <_{\cal L} v \leq_{\cal L} u$ \ then \
$sr_u <_{\cal L} B_{\cal R}(u, v)$. \\  
(3) \ If \ $su \ \frac{<}{>}\!\!\!|\,_{_{\cal L}} \ v$ \ then 
$sr_u \ \frac{<}{>}\!\!\!|\,_{_{\cal L}} \ B_{\cal R}(u, v)$. \\
(4) \ If \ $u \geq_{\cal L} v$ \ then
$B_{\cal R}(u, sv) = s \cdot B_{\cal R}(u, v)$. \\ 
(5) \ Analogous properties hold for $B_{\cal  L}$.
\end{Lm}
{\it Proof.} \ (1) By definition of $B_{\cal R}$ we have
$B_{\cal R}(u, v) = xr_u$ where $x$ is such that $v = xu$. But $v = asu$ 
for some $a$ since $v \leq_{\cal L} su$, hence we can pick $x = as$. So, \
$B_{\cal R}(u, v) = asr_u \leq_{\cal L} sr_u$.

\smallskip

By definition of $B_{\cal R}$ we have
$B_{\cal R}(sr_u, B_{\cal R}(u, v)) = xr_{sr_u}$, 
\  where $x$ is any element of $S$ such that \ $B_{\cal R}(u, v) = xsr_u$. 
 
Also, by definition of $B_{\cal R}$ we have 
$B_{\cal R}(su, v) = yr_{su}$, \ where $y$ is such that $v = ysu$.
By Lemma \ref{Lm1}, multiplying $v = ysu$ by $\alpha'$ we obtain \ 
$B_{\cal R}(u, v) = v\alpha' = ysu\alpha' = ysr_u$.
Thus, $B_{\cal R}(u, v) = ysr_u$, and since $x$ was any element such that
\ $B_{\cal R}(u, v) = xsr_u$, we can assume $x = y$. So,
\ $B_{\cal R}(su, v) = xr_{su}$. 
Moreover, \ $r_{sr_u} = r_{su}$ \ since $u \equiv_{\cal R} r_u$.
The result now follows.

\medskip

(2) \ By definition of $B_{\cal R}$ we have
\  $B_{\cal R}(u, v) = xr_u$ where $x$ is such that
$v = xu$. Hence $B_{\cal R}(u, v) = xr_u = v\alpha'$ \ where $\alpha'$ is
such that $u\alpha' = r_u$. Moreover, $v \geq_{\cal L} su$ 
(or $v >_{\cal L} su$),  thus
$B_{\cal R}(u, v) = v\alpha' \geq_{\cal L} su\alpha' = sr_u$ \ (or 
$>_{\cal L} su \alpha' = sr_u$).

\smallskip

By definition of $B_{\cal R}$ we have
$B_{\cal R}(B_{\cal R}(u, v), sr_u) = xr_{B_{\cal R}(u, v)}$, \ where $x$ is
such that $sr_u = x \cdot B_{\cal R}(u, v)$. \ 
By Lemma \ref{Lm1}, if we multiply the last equality by $\alpha$ we obtain
\ $su = xv$. 

By definition we also have \ $B_{\cal R}(v, su) = yr_v$, \  where $y$ is any 
element of $S$ such that $su = yv$. But we proved that $x$ also satisfies
\ $su = xv$. 
Thus we can assume $x = y$. 

So we have \ $B_{\cal R}(B_{\cal R}(u, v), sr_u) = yr_{B_{\cal R}(u, v)}$.
Moreover, since $B_{\cal R}(u, v) \equiv_{\cal R} v$ (by Lemma \ref{Lm3}),
we obtain the result.

\medskip

(3) \ This follows directly from Lemma \ref{Lm1}.

\medskip

(4) \  By definition, $B_{\cal R}(u, v) = yr_u$, where $yu = v$. \ Also
$B_{\cal R}(u, sv) = xr_u$, where $x$ is any element of $S$ such that 
$xu = sv$. Since $yu = v$, we have $syu = sv$, hence we can pick $x$ to be
$sy$. The result then follows.
\ \  $\Box$

\begin{Lm} \label{Lm8}
If \ $w \ \frac{<}{>}\!\!\!|\,_{_{\cal L}} \  s$ \ then 
$B(u, v, w) \ \frac{<}{>}\!\!\!|\,_{_{\cal L}} \  s$.  \ Similarly, \  
if \ $s \ \frac{<}{>}\!\!\!|\,_{_{\cal R}} \ u$ \ then 
$s \ \frac{<}{>}\!\!\!|\,_{_{\cal R}} \ B(u, v, w)$. 
\end{Lm} 
{\it Proof.} \ By contraposition, assume \  
\ $B(u, v, w) \ \frac{<}{>}\!\!\,_{_{\cal L}} \  s$. \ By definition, 
$B(u, v, w) = ux$, where $x$ is such that $w = vx$. Since $B(u, v, w)$ exists,
$u \leq_{_{\cal L}} v \geq_{_{\cal R}} w$; \ so $u = yv$ for some $y$.

Now we have $s \ \frac{<}{>}\!\!\,_{_{\cal L}} \ B(u, v, w) = ux = yvx =
yw \leq_{_{\cal L}} w$. \\
In case $s \leq_{_{\cal L}} ux$, the above implies \ $s \leq_{_{\cal L}} w$. \\
In case $s \geq_{_{\cal L}} ux$, the above implies \ $s \geq_{_{\cal L}} ux 
\leq_{_{\cal L}} w$, and hence, by unambiguity of the ${\cal L}$-order, \ 
$s \ \frac{<}{>}\!\!\,_{_{\cal L}} \ w$.

In either case, $s \ \frac{<}{>}\!\!\,_{_{\cal L}} \ w$.
\ \  $\Box$

\begin{Lm} \label{Lm9}
(Lemma 1.1.(5) in \cite{GrPaper}.) \ 
If \  $u \leq_{_{\cal L}} v \geq_{_{\cal R}} w \leq_{_{\cal L}} s 
\geq_{_{\cal R}} t$, \ then \  $B(u, v, w)  \leq_{_{\cal L}} s
\geq_{_{\cal R}} t$, \  $u \leq_{_{\cal L}} v \geq_{_{\cal R}} B(w, s, t)$, 
\  and \  $B(B(u, v, w),s, t) = B(u, v, B(w, s, t))$.
\end{Lm}
{\it Proof.} \ We have $B(u, v, w) \leq_{_{\cal L}} w$ \ by the definition 
of $B$, and
$w \leq_{_{\cal L}} s \geq_{_{\cal R}} t$, \ by assumption.
Also, \ $u \leq_{_{\cal L}} v \geq_{_{\cal R}} w$ \
by assumption, and \ $w \geq_{_{\cal R}} B(w, s, t)$ \ by Lemma \ref{Lm5}.
So the claimed order relations hold.

By Lemma \ref{Lm5}, $B(u, v, w) = yw$, where $u = yv$, \ and
by definition, \ $B(w, s, t) = wx$, where $t = sx$.    
Then by definition \ $B(u, v, B(w, s, t)) =  B(u, v, wx) = B(u, v, w) \cdot x$
(the latter equality holds by Lemma \ref{Lm6}).  
This is equal to $yw \cdot x$. 
A similar reasoning shows that $B(B(u, v, w),s, t)$ is also equal to $ywx$.
\ \ $\Box$

\begin{Lm} \label{Lm10} \ Assume that \
 \ $u \leq_{_{\cal L}} v \equiv_{_{\cal R}} w \geq_{_{\cal L}} s$ \ and \  
$c \in S - \{0\}$. Then: \\   
(a) \ \ \ \ \ \ $cs = B(u, v, w)$ \ iff \ $u = c \cdot B(s, w, v)$, \\  
(b) \ \ \ \ \ \ $cu = B(s, w, v)$ \ iff \ $s = c \cdot B(u, v, w)$.
\end{Lm}
{\it Proof of (a). \ } 
By Lemma \ref{Lm5}, there exist $x, y, x', y' \in S^{1}$ such that \\ 
$B(u, v, w) = ux = yw$, \ \  $w = vx$, \ $u = yv$ \ \  and \\ 
$B(s, w, v) = sx' = y'v$, \ \  $v = wx'$, \ $s = y'w$.

\smallskip
 
If the left side of the equivalence holds then \   
$yw = B(u, v, w) = cs = cy'w$,  so if we multiply by $x'$ we obtain 
 \ $u =  ywx' = cy'v = c \cdot B(s, w, v)$.

\smallskip

If the right side of the equivalence holds then \  
$u = c \cdot B(s, w, v) = cy'v$, \ so if we multiply by $x$ we obtain 
 \ $B(u, v, w) =  ux = cy'vx = cy'w = cs$.

\smallskip

The proof of {\it (b)} is similar. \ \ \ 
$\Box$

\begin{Lm} \label{Lm11} \ Assume that \ 
 \ $u \leq_{_{\cal L}} v \equiv_{_{\cal R}} w \geq_{_{\cal L}} s$. Then: 

\smallskip

\noindent (1) \ \ \ \ 
$B(u, v, w) \leq_{_{\cal L}} s$ \ iff \ $u \leq_{_{\cal L}} B(s, w, v)$. \  
The same holds with $\leq_{_{\cal L}}$ replaced by 
$>_{_{\cal L}}$ or \, $\frac{<}{>}\!\!\!|\,_{_{\cal L}}$.

\smallskip

\noindent (2. $\leq$) \ \ If \ $B(u, v, w) \leq_{_{\cal L}} s$ \ then \ 
$r_s = r_{B(s, w, v)}$ \ and \ $B_{_{\cal R}}(s, B(u, v, w)) = 
B_{_{\cal R}}(B(s, w, v), u)$. 

\noindent (2. $>$) \ \ If \ $B(u, v, w) >_{_{\cal L}} s$  \ then \ 
$r_u = r_{B(u, v, w)}$ \ and \ $B_{_{\cal R}}(B(u, v, w), s) = 
B_{_{\cal R}}(u, B(s, w, v))$.

\smallskip

Analogous properties hold for $B_{_{\cal L}}$.
\end{Lm}
{\it Proof.} (1): \  For $\leq_{_{\cal L}}$ this is an immediate consequence 
of Lemma \ref{Lm10} {\it (a)}. The result (1) for $\geq_{_{\cal L}}$ follows 
from Lemma \ref{Lm10} {\it (b)}. Since $>_{_{\cal L}}$ holds iff we have 
$\geq_{_{\cal L}}$ and not $\leq_{_{\cal L}}$,  we also obtain (1) for 
$>_{_{\cal L}}$. Also, since \ $\frac{<}{>}\!\!\!|\,_{_{\cal L}}$ \ holds 
iff we have neither $\leq_{_{\cal L}}$ nor $\geq_{_{\cal L}}$, we obtain
(1) for \ $\frac{<}{>}\!\!\!|\,_{_{\cal L}}$ .

\medskip

\noindent (2. $\leq$): \ If \ $B(u, v, w) \leq_{_{\cal L}} s$ \ then \
$B(s, w, v) = sx' \leq_{_{\cal R}} s$, and 
$B(s, w, v) \geq_{_{\cal R}} B(s, w, v) \cdot x = y'vx = y'w = s$, 
where $x', x, y'$ are as at the beginning of the proof of Lemma  \ref{Lm10}.
Thus \ $s \equiv_{_{\cal R}} B(s, w, v)$.

By definition, \ $B_{_{\cal R}}(s, B(u, v, w)) = x''r_s$, \ for any $x''$ 
such that \ $x''s = B(u, v, w)$. \\ 
And \ $B_{_{\cal R}}(B(s, w, v), u) = y''r_{B(s, w, v)} = y''r_s$, \ for any 
$y''$ such that \ $y'' \cdot B(s, w, v) = u$. \\ 
But by Lemma \ref{Lm10}, \ $x''s = B(u, v, w)$ \ iff \ 
$u = x'' \cdot B(s, w, v)$. \ So we can choose $y''$ to be $x''$. 
Then the equality follows.

\medskip

\noindent (2. $>$): \ The proof is very similar to that of (2. $\leq$). 
\ \ $\Box$

\begin{Lm} \label{Lm12} \ \ Assume that \
 \ $u \leq_{_{\cal L}} v \geq_{_{\cal R}} w \geq_{_{\cal L}} s$, \ and let
\ $c \in S$. \ Then:

\smallskip

\noindent (1) \ \ \ \ $B(u, v, w) = cs$ \ iff \ $B(u, v, r_w) = c \cdot 
B_{_{\cal R}}(w,s)$. \

\smallskip

\noindent (2) \ \ \ \ $c \cdot B(u, v, w) = s$ \ iff \ $c \cdot B(u, v, r_w) = 
B_{_{\cal R}}(w,s)$. \

\smallskip

Analogous properties hold for $B_{_{\cal L}}$.
\end{Lm}
{\it Proof.}  (1): \ Assume $B(u, v, w) = cs$, \ where
(by Lemma \ref{Lm5}) \ $B(u, v, w) = yw$ with $u = yv$. \ Multiplying
$yw = cs$ \ on the right by $\alpha'$, where $\alpha'$ is such
that $w\alpha' = r_w$, we obtain: \   
$yr_w = cs\alpha'$.

The left side $yr_w$ is equal to $B(u, v, r_w)$ by Lemma \ref{Lm5}, 
since $u = yv$. \ On the other hand, by the definition of 
$B_{_{\cal R}}$ we have, \ $B_{_{\cal R}}(w,s) = xr_w$
with $s = xw$.  Since $w\alpha' = r_w$, we have $B_{_{\cal R}}(w,s) =
xw\alpha' = s\alpha'$, which when multiplied by $c$ yields the right side.

\smallskip

Conversely, if \ $B(u, v, r_w) = c \cdot B_{_{\cal R}}(w,s)$ \ we will have
by Lemma \ref{Lm5} and by the definition of $B_{_{\cal R}}$, in the above 
notation: \ \  $yr_w = cs\alpha'$. 

Multiplying on the right by $\alpha$ (where $\alpha$ 
is such  that $r_w\alpha = w$), we obtain: \ $yw = cs\alpha'\alpha = cs$.  
We have \ $s\alpha'\alpha = s$ \ because we assumed \ $w >_{_{\cal L}} s$.
Thus \ $B(u, v, w) = yw = cs\alpha'\alpha = cs$. 

\smallskip

The proof of (2) is quite similar to the proof of (1). \ \  $\Box$

\begin{Lm} \label{Lm13} \ \ Assume that \
 \ $u \leq_{_{\cal L}} v \geq_{_{\cal R}} w \geq_{_{\cal L}} s$. \ Then:

\smallskip

\noindent (1) \ \ $B(u, v, w) \leq_{_{\cal L}} s$ \ iff \ 
$B(u, v, r_w) \leq_{_{\cal L}} B_{_{\cal R}}(w,s)$. \\   
The same is true with $\leq_{_{\cal L}}$ replaced by 
$>_{_{\cal L}}$  or \ $\frac{<}{>}\!\!\!|\,_{_{\cal L}}$ .

\smallskip

\noindent (2. $\leq$) \ \ If \ $B(u, v, w) \leq_{_{\cal L}} s$ \ then \
$s \equiv_{_{\cal R}} B_{_{\cal R}}(w,s)$ \ and  

$B_{_{\cal R}}(s, B(u,v,w)) = B_{_{\cal R}}(B_{_{\cal R}}(w,s), B(u,v,r_w))$.

\smallskip

\noindent (2. $>$) \ \ If \ $B(u, v, w) >_{_{\cal L}} s$  \ then \
$B(u, v, w) \equiv_{_{\cal R}} B(u, v, r_w)$ \ and \

$B_{_{\cal R}}(B(u, v, w), s) =
B_{_{\cal R}}(B(u, v, r_w), B_{_{\cal R}}(w,s))$.

\medskip

Analogous properties hold for $B_{_{\cal L}}$ : 

If \ $s \leq_{_{\cal R}} u \leq_{_{\cal L}} v \geq_{_{\cal R}} w$ \ then :

\smallskip

\noindent (1) \ \ $s \leq_{_{\cal R}} B(u, v, w)$ \ iff \ 
$B_{_{\cal L}}(s, u) \leq_{_{\cal R}} B(\ell_u, v, w)$. \\   
The same is true with $\leq_{_{\cal R}}$ replaced by
$>_{_{\cal R}}$  or \ $\frac{<}{>}\!\!\!|\,_{_{\cal R}}$ .

\smallskip

\noindent (2. $\leq$) \ \ If \ $s \leq_{_{\cal R}} B(u, v, w)$ \ then \  
$B(u, v, w) \equiv_{_{\cal L}} B(\ell_u, v, w)$ \ and \  

$B_{_{\cal L}}(s, B(u, v, w)) = 
B_{_{\cal L}}(B_{_{\cal L}}(s, u), B(\ell_u, v, w))$. 

\smallskip

\noindent (2. $>$) \ \ If \ $s >_{_{\cal R}} B(u, v, w)$ \ then \  
$s \equiv_{_{\cal L}} B_{_{\cal L}}(s, u)$  \ and \  

$B_{_{\cal L}}(B(u, v, w), s) = 
B_{_{\cal L}}(B(\ell_u, v, w), B_{_{\cal L}}(s, u))$.
\end{Lm}
{\it Proof.}  (1): \  The result for $\leq_{_{\cal L}}$ follows immediately
from Lemma \ref{Lm12} (1). From Lemma \ref{Lm12} (2), we have the corresponding
result for $\geq_{_{\cal L}}$. Combining the two we obtain the result for
$>_{_{\cal L}}$ \  and for \ $\frac{<}{>}\!\!\!|\,_{_{\cal L}}$ .

\medskip 

\noindent (2. $\leq$): \ By Lemma \ref{Lm3} we have \  
$r_s = r_{B_{_{\cal R}}(w,s)}$.  

\smallskip

We will apply Lemma \ref{Lm7} (2), which we quote here with different
parameters: 

{\it If} \ $s_ou_o \leq_{\cal L} v_o \leq_{\cal L} u_o$ \ {\it then}
$B_{\cal R}(v_o, s_ou_o) = B_{\cal R}(B_{\cal R}(u_o, v_o), s_or_{u_o})$.  

\noindent Let \ $v_o = s$, \ $u_o =  y$, \ and \ $s_o = w$, \ where 
(by Lemma \ref{Lm5}),
\ $B(u, v, w) = yw$ \ and \ $B(u, v, r_w) = yr_w$ \ with $yv = u$. Then
\ $s_ou_o = B(u, v, w)$ \ and \ $s_or_{u_o} = B(u, v, r_w)$. \ By assumption,
$B(u, v, w) \leq_{_{\cal L}} s <_{_{\cal L}} w$, so 
$s_ou_o \leq_{\cal L} u_o \leq_{\cal L} v_o$, 
\ hence  Lemma \ref{Lm7} (2) is indeed applicable here.  
By substituting, the claimed result then follows immediately.

\medskip

\noindent (2. $>$): \ By Lemma \ref{Lm5} we have \ $B(u, v, w) = yw$ \ and 
\ $B(u, v, r_w) = yr_w$, \ with $u = yv$. \ Since $w \equiv_{_{\cal R}} r_w$
we obtain \ $B(u, v, w) \equiv_{_{\cal R}} B(u, v, r_w)$.

\smallskip

We will apply Lemma \ref{Lm7} (1), which we quote here with different
parameters:

{\it If} \ $u_o \geq_{\cal L} s_ou_o \geq_{\cal L} v_o$ \ {\it then}
$B_{\cal R}(s_ou_o, v_o) = B_{\cal R}(s_or_{u_o}, B_{\cal R}(u_o, v_o))$.

\noindent Let \  $s_o = y$, \ and \  $u_o = w$, \ where $B(u, v, w) = yw$ 
and $B(u, v, r_w) = yr_w$, with $yv = u$ (by Lemma \ref{Lm5}). And let 
\ $v_o = s$. \ Since by our assumptions \  $w >_{\cal L} B(u, v, w) >_{\cal L} 
s$, Lemma \ref{Lm7} (1) can be applied. The claimed result then follows 
immediately by substitution. \ \  $\Box$

\begin{Lm} \label{Lm14} \ \ Assume that \
 \ $u \leq_{_{\cal L}} v \geq_{_{\cal R}} w$. \  Then \  
$B_{\cal R}(v,u) \leq_{_{\cal L}} r_v \geq_{_{\cal R}} w$ \ and \  
$B(B_{\cal R}(v,u), r_v, w) = B(u, v, w)$. 

\smallskip

Analogous properties hold for $B_{_{\cal L}}$:

If \ $u \leq_{_{\cal L}} v \geq_{_{\cal R}} w$ \ then \   
$u \leq_{_{\cal L}} \ell_v \geq_{_{\cal R}} B_{_{\cal L}}(w, v)$ \ and \  
$B(u, v, w) = B(u, \ell_v, B_{_{\cal L}}(w, v))$.
\end{Lm}
{\it Proof.} \ The fact that $B_{\cal R}(v,u) \leq_{_{\cal L}} r_v 
\geq_{_{\cal R}} w$ \ is obvious from the definition of $B_{\cal R}$. 

\  By Lemma \ref{Lm5}, \ $B(u, v, w) = x_1w$ \ for any $x_1$ such
that $u = x_1v$. \ Also, by definition, \ $B_{\cal R}(v,u) = x_2r_v$ \ for any
$x_2$ such that $u = x_2v$; \ therefore we can choose $x_2 = x_1$.
 
Now $B(B_{\cal R}(v,u), r_v, w) = B_{\cal R}(v,u) \, z$ \ with $w = r_vz$, \  
hence \ $B(B_{\cal R}(v,u), r_v, w)  = x_1r_vz = x_1w$. \ This proves the 
result. \ \  $\Box$

\begin{Lm} \label{Lm15} \ \ Assume that \
 \ $u \leq_{_{\cal L}} v \geq_{_{\cal R}} w \leq_{_{\cal L}} s$. \ Then \ 
$B_{\cal R}(s, B(u, v, w)) = B(u, v, B_{\cal R}(s,w))$. 

\smallskip

Analogous properties hold for $B_{_{\cal L}}$: 

If \ $s \geq_{_{\cal R}} u \leq_{_{\cal L}} v \geq_{_{\cal R}} w$ \ then \  
$B_{_{\cal L}}(B(u, v, w), s) = B(B_{_{\cal L}}(u,s), v, w)$. 
\end{Lm}
{\it Proof.} \  By definition, \ $B_{\cal R}(s, B(u, v, w)) = x_1r_s$ \ 
where \ $x_1s = B(u, v, w) = uz$, with (by definition of $B$) $w = vz$.
We also have: \\  
$B(u, v, B_{\cal R}(s,w))$ \\ 
$= y \, B_{\cal R}(s,w)$ \hspace{0.8in}  where $y$ is such that \ $u = yv$ \\
$= yx_2r_s$ \ \hspace{1.in} where \ $x_2$ is such that $x_2s = w$ \\ 
$= yx_2s\alpha'$ \ \hspace{1.in} where \ $\alpha'$ is such that $r_s = s\alpha'$\\ 
$= yw\alpha'$ \ \hspace{1.in} since \ $x_2s = w$ \\ 
$= yvz\alpha'$ \ \hspace{1.in} since \ $w = vz$ \\  
$= uz\alpha'$ \ \hspace{1.in}  since \ $u = yv$ \\
$= B(u, v, w)\alpha'$ \\ 
$= x_1s\alpha'$ \\  
$= x_1r_s$ \\  
$= B_{\cal R}(s, B(u, v, w))$ \hspace{0.7in} as we saw in the beginning of 
this proof.           \ \ \ \ $\Box$

\begin{Lm} \label{Lm16} \ \ Assume that \
 \ $u \geq_{_{\cal L}} v \geq_{_{\cal R}} w$. \ Then   \\ 
(1) \ \ \ \ 
$B_{_{\cal R}}(u, v) \equiv_{_{\cal L}} B_{_{\cal R}}(u, \ell_v))$, \\  
(2) \ \ \ \ $B_{_{\cal L}}(w, B_{_{\cal R}}(u, v)) = 
B_{_{\cal L}}(B_{_{\cal L}}(w,u), B_{_{\cal R}}(u, \ell_v))$. 
\end{Lm}
{\it Proof.} \ Property (1) follows easily from Lemma \ref{Lm7} (4). 

\smallskip

\noindent (2): \ \ \ \   Let $\beta$ and $\beta'$ be such that 
$v = \beta \ell_v$ and $\ell_v = \beta' v$. \ 
By definition, \ $B_{_{\cal R}}(u, v) = xr_u$, where
$xu = v$. Hence, by the definition of $B_{_{\cal R}}$, we have
\ $B_{_{\cal R}}(u, \ell_v) = \beta' xr_u$ \ since $\beta' x$ 
satisfies $\beta' xu = \ell_v$. 

Thus, \ $B_{_{\cal L}}(w, B_{_{\cal R}}(u, v)) = B_{_{\cal L}}(w, xr_u) = 
\ell_{xr_u}y_1$, \ where $y_1$ is such that $w = xr_uy_1$. 

On the other hand, \ 
$B_{_{\cal L}}(B_{_{\cal L}}(w,u), B_{_{\cal R}}(u, \ell_v)) = 
\ell_{B_{_{\cal R}}(u, \ell_v)} y_2 = \ell_{xr_u}y_2$, \ since 
$B_{_{\cal R}}(u, v) \equiv_{_{\cal L}} B_{_{\cal R}}(u, \ell_v))$ \  
(as we just proved in (1)).
\ Here, by the definition of $B_{_{\cal L}}$, \ 
$y_2$ is any element of $S$ such that \ 
$B_{_{\cal L}}(w,u) = B_{_{\cal R}}(u, \ell_v) \, y_2$.
We saw that the latter is equal to \ $\beta' xr_u y_2$.
By the definition of $B_{_{\cal L}}$ we also have \ $B_{_{\cal L}}(w,u) = 
\ell_vy_3$ where $y_3$ is such that $w =vy_3$.

Therefore \ $\ell_vy_3 = \beta' xr_u y_2$. \ Multiplying on the left by 
$\beta$ yields \ $w = vy_3 = xr_u y_2$, \ i.e., \ $y_2$ satisfies 
\ $w = xr_u y_2$, which is the defining property of $y_1$. 

Hence, $y_2$ can be chosen above so that $y_2 = y_1$.
\ \ \ \ $\Box$

\begin{Lm} \label{Lm17} \ \ Assume that \
 \ $u' \geq_{_{\cal L}} v \leq_{_{\cal R}} w$. \ Then   \ \ \
$B_{_{\cal L}}(B_{_{\cal R}}(u, v), w) = B_{_{\cal R}}(u, B_{_{\cal L}}(v,w))$.
\end{Lm}
{\it Proof.} \ By the definition of $B_{_{\cal R}}$ and $B_{_{\cal L}}$, \  
$B_{_{\cal R}}(u, v) = xr_u$, \ where $v = xu$, \ and \ 
$B_{_{\cal L}}(v,w) = \ell_wy$, \ where $v = wy$. \ Let $\alpha$, $\alpha'$,
$\beta$ and $\beta'$
be such that $r_u \alpha = u$, $u\alpha' = r_u$, 
$\beta \ell_w = w$, and $\beta' w = \ell_w$.  

Then
$B_{_{\cal L}}(B_{_{\cal R}}(u, v), w) = \ell_wy_1$, \ where $y_1$ is such 
that \ $(xr_u =) \ B_{_{\cal R}}(u, v) = wy_1$. 

Also, \ $B_{_{\cal R}}(u, B_{_{\cal L}}(v,w)) = x_1r_u$, \ where $x_1$ is such
that \ $(\ell_wy =) \ B_{_{\cal L}}(v,w) = x_1u$. \ By multiplying 
the latter equalities by $\beta$ we obtain: \ 

\noindent (*) \hspace{1.2in} $wy = \beta x_1u$.

\smallskip

\noindent We need to show that \ \ $\ell_wy_1 =x_1r_u$. 

We saw that \ $v = xu = xr_u \alpha  = B_{_{\cal R}}(u, v) \, \alpha$ \ 
(by the choice of $x$ and of $\alpha$, and by the definition of 
$B_{_{\cal R}}$). Thus \ 

\hspace{1.2in} \ \ \ $B_{_{\cal R}}(u, v) \, \alpha = v$. 

\noindent In this equation we replace $v$ by  $wy$ (see the definition of 
$B_{_{\cal L}}(v,w)$), and we replace $B_{_{\cal R}}(u, v)$ by $wy_1$
(see the expression for $B_{_{\cal L}}(B_{_{\cal R}}(u, v), w)$).
Thus, 

\hspace{1.2in} $wy_1\alpha = wy$.

\smallskip

\noindent By (*) we can replace $wy$ by $\beta x_1u$.  \ So,

\hspace{1.2in} \ \ $wy_1\alpha = \beta x_1u$.

\noindent Multiplying this by $\alpha'$ (on the left) and by 
$\beta'$ (on the right) yields \ 
$\ell_wy_1 = x_1r_u$, \ which is what we wanted.
\ \ \ $\Box$ 

\begin{Lm} \label{Lm18} \ \ Assume that \ \ 
$u \leq_{_{\cal L}} v \geq_{_{\cal R}} w$. \ Then   \ \ \
$B(B_{_{\cal R}}(v,u), r_v, w) = B(u, \ell_v, B_{_{\cal L}}(w,v))$.
\end{Lm}
{\it Proof.} \  By the definition of $B_{_{\cal R}}$ and $B_{_{\cal L}}$, we 
have: \\ 
$B(B_{_{\cal R}}(v,u), r_v, w) = B(xr_v, r_v, w)$, \ where $u = xv$, and  \\ 
$B(u, \ell_v, B_{_{\cal L}}(w,v)) = B(u, \ell_v, \ell_v y)$, \ where $w = vy$.

\smallskip

By the definition of $B$, \ $B(xr_v, r_v, w) = xr_vz_1$, \ 
where $w = r_vz_1$.  \ Hence, \ $B(xr_v, r_v, w) = xw$.  

Similarly, \ $B(u, \ell_v, \ell_v y) = uz_2$, \ where $z_2$ is any element of 
$S$ satisfying $\ell_v y = \ell_v z_2$; hence we can pick $z_2$ to be $y$. 
\ Then we have \ $B(u, \ell_v, \ell_v y) = uy = xvy$ \ (since $u = xv$), \ 
and \ $xvy = xw$ \ (since $vy = w$). Thus \ $B(u, \ell_v, \ell_v y) = xw$, \ 
which is equal to  \ $B(xr_v, r_v, w)$, \ as we saw. \ \ \ $\Box$

\section{Termination}

In this section we prove that the rewrite system for $(S)_{\mathrm{reg}}$ is 
terminating.

\begin{Lm}
If the sub-system consisting of the rules (2.1)--(2.4) is terminating then
the whole rewrite system is terminating.
\end{Lm}
{\it Proof.} \  Imagine, by contraposition, that the whole rewrite system 
allows an infinite rewrite chain. Since the first group of rules is strictly 
length-reducing, the chain contains only rules of the form (2.1)--(2.4), 
from some point on. Hence the rules (2.1)--(2.4) do not form a 
terminating system.   \ \ \ $\Box$

\medskip

The rest of this section deals with the proof that the
sub-system consisting of the rules (2.1)--(2.4) terminates. In the remainder of 
this section, rewriting means applying the rules (2.1)--(2.4).  

\medskip

Since the rules (2.1)--(2.4) are length-preserving, the notion of 
{\em position} in a string is invariant under rewriting. 
More precisely, a string $x = (x_1, \ldots,x_n)$  of length $n$
over the generators of $(S)_{\mathrm{reg}}$ has positions $1, 2, \ldots,
n$, and when a rule of type (2.1)--(2.4) is applied, the new string still 
has positions $1, 2, \ldots, n$. 

Our first step is to find factorizations of strings that are preserved under 
rewriting. See \cite{Bi3} for more background on preserved factorization schemes;
here we do not need exact definitions since the context will make everything clear.

\begin{Lm}
In a string, a position occupied by 0 is invariant under rewriting. 
Similarly, the fact that a position is occupied by an element of $S - \{0\}$ 
(respectively by an element of $\overline{S - \{0\}}$) is invariant 
under rewriting.
\end{Lm}
{\it Proof.} \ 
Since the rules (2.1)--(2.4) do not use the symbol 0, a position 
occupied by 0 will never change, and a non-0 symbol never turns into 0. 
Similarly, a position occupied by an element $s \in S- \{0\}$ will always 
remain occupied by an element of $S- \{0\}$, although the value of $s$ can 
change. Similarly for $\overline{S - \{0\}}$.
\ \ $\Box$

\begin{Lm} (Preservation of $<_{_{\cal L}}$, $\equiv_{_{\cal L}}$, 
$>_{_{\cal L}}$, and $\frac{<}{>}\!\!\!|\,_{_{\cal L}} \ $, and similarly for 
${\cal R}$). \\ 
In a string, a pair of positions occupied by elements \ 
$(s, \overline{t}) \in S \times \overline{S}$ \ with $s <_{_{\cal L}} t$ 
(or $\equiv_{_{\cal L}}$ or $>_{_{\cal L}}$ or 
$\frac{<}{>}\!\!\!|\,_{_{\cal L}}$) 
will always remain occupied by some pair of 
in $S \times \overline{S}$ related by $<_{_{\cal L}}$ (respectively 
$\equiv_{_{\cal L}}$ or $>_{_{\cal L}}$ or $\frac{<}{>}\!\!\!|\,_{_{\cal L}}$). 
Similarly, for a pair in $\overline{S} \times S$ related by 
$<_{_{\cal R}}$ (or $\equiv_{_{\cal R}}$ or $>_{_{\cal R}}$ or 
$\frac{<}{>}\!\!\!|\,_{_{\cal R}}$), this relation is preserved between these
two positions. 
\end{Lm}
{\it Proof.} \  
Let us look at the four ways $s$ or $\overline{t}$ could be changed
when a rule is applied just to the left or right of $(s, \overline{t})$.

If the symbol to the left of $(s, \overline{t})$ is $\overline{u}$, with
$u >_{_{\cal R}} s$, then (2.2) can change $(\overline{u}, s, \overline{t})$
into $(\overline{\ell_u}, B_{_{\cal L}}(s, u), \overline{t})$. 
Since $B_{_{\cal L}}(s, u) \equiv_{_{\cal L}} s$ (by Lemma \ref{Lm3}), we still have 
$B_{_{\cal L}}(s, u) <_{_{\cal L}} t$ at this pair of positions.

If the symbol to the left of $(s, \overline{t})$ is $\overline{u}$, with
$u \leq_{_{\cal R}} s$, then (2.3)  can change $(\overline{u}, s, \overline{t})$
into $(\overline{B_{_{\cal L}}(u, s)}, \ell_s, \overline{t})$. 
Since $\ell_s \equiv_{_{\cal L}} s$ we still have $\ell_s <_{_{\cal L}} t$
at this pair of positions.

If the symbol to the right of $(s, \overline{t})$ is $v$ with $t >_{_{\cal R}} v$
(or $t \leq_{_{\cal R}} v$) then the reasoning is similar.
\ \ $\Box$

\bigskip

As a consequence of these preservation lemmas we can factor any string
into {\it maximal subsegments}, defined by the following properties: \\
$\bullet$ \ \ 0 does not occur in a subsegment, unless the subsegment consists
of only 0; \\
$\bullet$ \ \ neighboring positions in a subsegment are occupied by pairs in 
$S \times \overline{S}$ or $\overline{S} \times S$; \\
$\bullet$ \ \ the incomparability relation $\frac{<}{>}\!\!\!|$ (for ${\cal L}$ 
or ${\cal R}$) does not occur inside a subsegment.\\
We call such subsegments {\em continuous strings}, i.e., we view the break 
between two maximal such subsegments as a discontinuity.
The rewrite rules (2.1)--(2.4) preserve this factorization; no rewrite rule
applies to two positions that are in different maximal subsegments.
 
A string is called continuous iff it consists of just one maximal subsegment.
For a continuous string $x = (x_1, \ldots, x_n)$ over the generators of
$(S)_{\mathrm{reg}}$ and a position $i$ ($1 \leq i < n$), we write $x_i > x_{i+1}$
(or $<$, $\leq$,$\geq$) iff the corresponding $\cal R$- or
$\cal L$-relation holds in $S$ according to the above Lemma. 

\bigskip

{\bf Definition.} Let $x = (x_1, \ldots, x_n)$ be a continuous string of length 
$n$. We call a position $i$ ($1 \leq i \leq n$) in $x$ {\em maximal} \ iff \\
$\bullet$ \ $i = 1$ and $x_1 > x_2$, \ or \\
$\bullet$ \ $i = n$ and  $x_{n-1} \leq x_n$, \ or \\
$\bullet$ \ $1 < i < n$ and  $x_{i-1} \leq x_i > x_{i+1}$. 

\medskip

By Lemma 4.3, maximal positions remain maximal during rewriting. 

\begin{Lm} (Maximal positions).  \\
During the rewriting of a continuous string using rules (2.1)--(2.4),
an element of $S \cup \overline{S}$ at a maximal position is rewritten at 
most twice. From then on, the symbol at the maximal position never changes. 
\end{Lm}
{\it Proof.} \ 
Suppose that a maximal position is occupied by an element $s \in S$ (the case
of an element of $\overline{S}$ is similar).
Let $\overline{u},s, \overline{v}$ be the neighboring elements in the 
continuous string, with $u \leq_{_{\cal R}} s >_{_{\cal L}} v$. The element
$\overline{u}$ or the element $\overline{v}$ may be absent. 
If (2.3) is applied, $(\overline{u}, s)$ will be rewritten to
$(\ldots , \ell_s)$. If (2.1) is applied, $(s, \overline{v})$ will be rewritten
to $(r_s, \ldots)$. If (2.3) is now applied (or (2.1) is applied to the 
previous alternative), the element at the maximal position is rewritten to 
$\ell_{r_s}$ (respectively $r_{\ell_s}$). Further rewriting with rules
(2.1), (2.3) cannot change the element at the maximal position because 
$r_{\ell_{r_s}} = \ell_{r_s}$ and $\ell_{r_{\ell_s}} = r_{\ell_s}$. This
follows from the special choice of the representatives
of the ${\cal L}$- and ${\cal R}$-classes; recall that 
$\equiv_{\cal H}$-related representatives are equal. 
\ \  $\Box$ 

Note that the above Lemma (and the termination property itself) is not true
if the representatives of the ${\cal L}$- and ${\cal R}$-classes are chosen 
differently than we did (except in trivial cases, e.g., when 
$S - \{0\}$ has no strict $>_{_{\cal R}}$ and $>_{_{\cal L}}$ chains).   

\begin{Lm} (Chains \ $\ldots > \cdot > \ldots$ \ and chains \ 
$\ldots \leq \cdot \leq \ldots$ \ stabilize). \\ 
If $s \in S$ occurs in a continuous string,
with $\ldots >_{_{\cal L}} s >_{_{\cal R}} \ldots$ \ or \ 
$\ldots \leq_{_{\cal R}} s \leq_{_{\cal L}} \ldots$ in this string,
then after a finite number of applications of the rules (2.1)--(2.4) to 
the string, the symbol at the position of $s$ will not change any more. 

The same is true for an occurrence of $\overline{s} \in \overline{S}$ 
in a continuous string, with 
$\ldots >_{_{\cal R}} s >_{_{\cal L}} \ldots$ \ or \   
$\ldots \leq_{_{\cal L}} s \leq_{_{\cal R}}  \ldots$.
\end{Lm}
{\it Proof.} \  Let us consider a continuous string $(\ldots, s, \ldots)$ 
with $s \in S$ and 
$\ldots >_{_{\cal L}} s >_{_{\cal R}} \ldots$. 
By the previous lemma, we know that
the element at the maximal position towards the left of $s$ will eventually
stabilize. By induction, suppose that all elements in the descending 
alternatining $>_{_{\cal L}}$--$>_{_{\cal R}}$ chain to the left of $s$ 
have stabilized.
No rule among (2.1)--(2.4) can be applied to the left of $s$ in this chain 
anymore (otherwise the
element just left of $s$ would change again, since $u \neq r_u$, resp. 
$u \neq \ell_u$ in the rules). On the other hand, 
if a rule is applied to $s$ and the element
just right of $s$ (in that case it would be rule (2.2)), then $s$ is replaced 
by $r_s$ and after this, no rule can be applied anymore at this position.

Let us also consider the case of a continuous string $(\ldots, s, \ldots)$ 
with $s \in S$ and 
$\ldots \leq_{_{\cal R}} s \leq_{_{\cal L}} \ldots$.  As before, let us assume
that all maximal positions have stabilized, and let us assume by
induction that all elements in the ascending alternatining 
$\leq_{_{\cal L}}$--$\leq_{_{\cal R}}$ chain to the right of $s$ have stabilized. 
Again, no rule will
be applied to the right of $s$ anymore. On the other hand, if  a rule is 
applied to $s$ and the element just left of $s$ (in that case it will be 
rule (2.3), then $s$ is replaced by $\ell_s$, and after this, no rule can be
applied anymore at this position.

The reasoning is similar in the other cases.  \ \  $\Box$  

\bigskip

{\bf Definition.} \  Let $x = (x_1, \ldots, x_n)$ be a continuous string of 
length $n$. We call a position $i$ (1 $\leq i \leq n$) 
{\em minimal} iff \\
$\bullet$ \ $i = 1$ and $x_1 \leq x_2$, \ or  \\
$\bullet$ \ $i = n$ and $x_{n-1} > x_n$, \ or \\
$\bullet$ \ $1 < i <n$ and $x_{i-1}> x_i \leq x_{i+1}$.

\medskip

By Lemma 4.3, minimal positions remain minimal during rewriting.

\begin{Lm} (Minimal positions stabilize). \\
After a finite number of applications of the rules (2.1)--(2.4) to a 
continuous string the symbols at the minimal positions do not change anymore.
\end{Lm}
{\it Proof.} \   Consider the case of a minimal position occupied by an 
element $v \in S - \{0\}$, occurring in a context 
$(\ldots, \overline{u}, v, \overline{w}, \ldots)$,  with $u >_{_{\cal R}} v 
\leq_{_{\cal L}} w$.
By the previous Lemma we assume that $u$ and $w$ will not change anymore. 
Then no rule can be applied to $v$, otherwise $u$ or $w$  would change again, since 
$s \neq r_s$, resp. $s \neq \ell_s$ in the rules. \ \ $\Box$

\bigskip

The Lemmas imply that all positions in a string eventually stabilize for the 
rewrite rules (2.1)--(2.4).

\section{Local confluence}

This section contains the proof that the rewrite system for $(S)_{\mathrm{reg}}$ is
locally confluent. We have to look at all the overlap cases (see \cite{Ja}),
which is tedious but straightforward in each case. Each case is either trivial
or it is resolved by using the properties of $B$, $B_{\cal L}$ and 
$B_{\cal R}$ proved in Section 3.

\bigskip

\noindent {\bf Overlap 1.1--1.1: \ }
$(s t, u) \ \stackrel{1.1}{\longleftarrow} \ (s, t, u) \ 
\stackrel{1.1}{\longrightarrow} \ (s, t u)$.

Then $(s t, u)\ \stackrel{1.1}{\longrightarrow} \ (s t u) \ 
\stackrel{1.1}{\longleftarrow} \
(s, t u)$, where we also use associativity
of the multiplication in $S$. 

\medskip

The overlap for the $\overline{S}$-form of rule 1.1 has the form \\
$(\overline{ts}, \overline{u}) \ \stackrel{1.1}{\longleftarrow} \ 
(\overline{s}, \overline{t}, \overline{u}) \ \stackrel{1.1}{\longrightarrow} \ 
\overline{s}, \overline{ut})$. \\
Confluence follows easily as above. 

\bigskip 

\noindent {\bf Overlaps with 1.2: \ } 
In all overlaps with rule 1.2 one easily 
shows confluence to (0).

\bigskip 

\noindent {\bf Overlap 1.1--1.3: \ }

{\bf Case 1.} \ $S$-form of rule 1.1.

\noindent 
$(t u, \overline{v})  \ \stackrel{1.1}{\longleftarrow} \ (t, u, \overline{v})$
$ \ \stackrel{1.3}{\longrightarrow} \  (t, 0)$  \ \ \ \ \  
where \  $u \ \frac{<}{>}\!\!\!|\,_{\cal L}\  v$.

\noindent Then \ $(t, 0) \ \stackrel{1.2}{\longrightarrow} \ (0)
\ \stackrel{1.3}{\longleftarrow} \ (tu, \overline{v})$. 
\ The last application of rule 1.3 is justified by the following.

\medskip 

\noindent {\it Claim: \ If \  $u \ \frac{<}{>}\!\!\!|\,_{\cal L}\  v$ \ 
then \ $tu \ \frac{<}{>}\!\!\!|\,_{\cal L}\  v$.}

\smallskip

\noindent Proof of the Claim: \  By contraposition, if 
\ $u \geq_{\cal L} tu \geq_{\cal L} v$ \ then obviously \ 
$u \geq_{\cal L}  v$. \ And if \ $u \geq_{\cal L} tu \leq_{\cal L} v$ \ 
then \ $u \ \frac{<}{>}_{\cal L}\  v$, \ by unambiguity of $S$. \ 
This proves the Claim.

\medskip

{\bf Case 2.} \ $\overline{S}$-form of rule 1.1.

\noindent 
$(0, \overline{v})  \ \stackrel{1.3}{\longleftarrow} \ 
(t, \overline{u}, \overline{v}) \ \stackrel{1.1}{\longrightarrow} \  
(t, \overline{vu})$, \ \ \ \  
where \ $t \ \frac{<}{>}\!\!\!|\,_{\cal L}\ u$. 

Confluence is proved in the same way as above.

\bigskip

\noindent {\bf Overlap 1.1--1.4: \ } Similar to the previous case.

\bigskip

\noindent {\bf Overlap 1.1--1.5:}  

{\bf Case 1.} \ \  $(tu, \overline{v}, w) \ 
\stackrel{1.1}{\longleftarrow} \ (t, u, \overline{v}, w) \ 
\stackrel{1.5}{\longrightarrow} \ (t, B(u, v, w))$, \ \ \  
where \ $u \leq_{\cal L} v \geq_{\cal R} w$.

\smallskip

Then \ $(tu, \overline{v}, w) \ \stackrel{1.5}{\longrightarrow} \ 
B(tu, v, w)$, and \  
$t \cdot B(u, v, w) \ \stackrel{1.1}{\longleftarrow} \ (t, B(u, v, w))$. \\
But by Lemma \ref{Lm6}, \ $B(tu, v, w) = t \cdot B(u, v, w)$, so we have 
confluence. 

\medskip

{\bf Case 2.} \ \ $(u, \overline{v}, w t) \ 
\stackrel{1.1}{\longleftarrow} \ (u, \overline{v}, w, t) \  
\stackrel{1.5}{\longrightarrow} \ (B(u, v, w) \cdot t)$ \ \ \  
where \ $u \leq_{\cal L} v \geq_{\cal R} w$.
 
\smallskip

As in the previous case, we have confluence by Lemma \ref{Lm6}.

Here we only considered the $S$-form of rule 1.1; the $\overline{S}$-form
does not overlap with 1.5.

\bigskip

\noindent {\bf Overlap 1.1--1.6: \ } Only the $\overline{S}$-form of 1.1 
overlaps with 1.6. Confluence is proved in a similar  way as in 
1.1--1.5.

\bigskip

\noindent {\bf Overlap 1.1($S$-form) -- 2.1: \ } \ $(su, \overline{v}) \ 
\stackrel{1.1}{\longleftarrow} \ 
(s, u, \overline{v}) \  \stackrel{2.1}{\longrightarrow} \ 
(s, r_u, \overline{B_{\cal R}(u, v)})$, \ \ \  
where \ $u >_{\cal L} v$.

\smallskip

{\bf Case 1. \ } \ $su  >_{\cal L} v$. \

Then \ $(su, \overline{v}) \ \stackrel{2.1}{\longrightarrow} 
(r_{su}, \overline{B_{\cal R}(su, v)})$, \ since \ $su  >_{\cal L} v$. \\
Moreover, 
$(s, r_u, \overline{B_{\cal R}(u, v)}) \ \stackrel{1.1}{\longrightarrow} 
(sr_u, \overline{B_{\cal R}(u, v)})  \ \stackrel{2.1}{\longrightarrow} 
(r_{sr_u}, \overline{B_{\cal R}(sr_u, B_{\cal R}(u, v))})$, \ 
where the latter application of rule 2.1 is justified since
$sr_u >_{\cal L} B_{\cal R}(u, v)$ \ (indeed we assumed $su  >_{\cal L} v$, 
so by Lemma \ref{Lm1}, \ 
$sr_u = su\alpha' >_{\cal L} v\alpha' = B_{\cal R}(u, v)$). 

\smallskip

To have confluence we need $r_{su} = r_{sr_u}$ (which easily follows from 
$u \equiv_{\cal R} r_u$), and \ 
$B_{\cal R}(su, v) = B_{\cal R}(sr_u, B_{\cal R}(u, v))$ \ (which is proved in 
Lemma \ref{Lm7} (1)). 

\medskip

{\bf Case 2. \ } \ $su \leq_{\cal L} v$. \ 

Then \ $(su, \overline{v}) \ \stackrel{2.4}{\longrightarrow}
(B_{\cal R}(v, su), \overline{r_v})$. 

Moreover, $(s, r_u, \overline{B_{\cal R}(u, v)}) \ 
\stackrel{1.1}{\longrightarrow} \  (sr_u, \overline{B_{\cal R}(u, v)}) \ 
\stackrel{2.4}{\longrightarrow} \  
(B_{\cal R}(B_{\cal R}(u, v),sr_u), \overline{r_{B_{\cal R}(u, v)}})$. \  
The latter application of rule 2.4 is justified since \ 
$sr_u \leq_{\cal L} B_{\cal R}(u, v)$, which follows from the assumption
$su \leq_{\cal L} v$ and from Lemma \ref{Lm1}. 

In order to have confluence we need \ $B_{\cal R}(B_{\cal R}(u, v),sr_u) =
B_{\cal R}(v, su)$ \ (which was proved in Lemma \ref{Lm7} (2)), and \  
$r_{B_{\cal R}(u, v)} = r_v$ \ (which follows from Lemma \ref{Lm3}). 

\medskip

{\bf Case 3. \ } \ $su \ \frac{<}{>}\!\!\!|\,_{_{\cal L}} \ v$.

Then \ $(su, \overline{v}) \ \stackrel{1.3}{\longrightarrow} \ (0)$. 

Moreover, $(s, r_u, \overline{B_{\cal R}(u, v)}) \
\stackrel{1.1}{\longrightarrow} \  (sr_u, \overline{B_{\cal R}(u, v)})$.
\ By Lemma \ref{Lm7} (3), \ $sr_u \ \frac{<}{>}\!\!\!|\,_{_{\cal L}} \ 
B_{\cal R}(u, v)$, \ so we can now apply rule 1.3, thus obtaining confluence 
to (0).

\medskip

\noindent {\bf Overlap 1.1($\overline{S}$-form) -- 2.1: \ } \ 
$(r_u, \overline{B_{\cal R}(u, v)}, \overline{s}) 
\ \stackrel{2.1}{\longleftarrow} \ (u, \overline{v}, \overline{s}) \ 
\stackrel{1.1}{\longrightarrow} \ (u, \overline{sv})$, \ \ \  
where \ $u >_{\cal L} v$.

\smallskip

Then \ $(r_u, \overline{B_{\cal R}(u, v)}, \overline{s}) \  
\stackrel{1.1}{\longrightarrow} \ (r_u, \overline{s \, B_{\cal R}(u, v)})$, \ 
and \ $(u, \overline{sv}) \ \stackrel{2.1}{\longrightarrow} \ 
(r_u, \overline{B_{\cal R}(u, sv)})$; 2.1 was applicable since 
$u >_{\cal L} v \geq_{\cal L} sv$. \ Confluence than follows directly from 
Lemma \ref{Lm7} (4).  

\bigskip

\noindent {\bf Overlap 1.1--2.2: \ } This is similar to the overlap 
1.1--2.1. 

\bigskip

\noindent {\bf Overlap 1.1--2.3: \ } This is similar to the overlap
1.1--2.4, which we consider next.

\bigskip

\noindent {\bf Overlap 1.1($S$-form) -- 2.4: \ } $(sv, \overline{u}) \ 
\stackrel{1.1}{\longleftarrow} \ (s, v, \overline{u}) \ 
\stackrel{2.4}{\longrightarrow} \ (s, B_{\cal R}(u, v), \overline{r_u})$, \ \ 
\   
where \ $v \leq_{\cal L} u$. 

Then \ $(sv, \overline{u}) \ \stackrel{2.4}{\longrightarrow} \ 
(B_{\cal R}(u, sv), \overline{r_u})$. 

Moreover, \ $(s, B_{\cal R}(u, v), \overline{r_u}) \  
\stackrel{1.1}{\longrightarrow} \ (s \cdot B_{\cal R}(u, v), \overline{r_u})$.

Confluence then follows from Lemma \ref{Lm7} (4).

\bigskip

\noindent {\bf Overlap 1.1($\overline{S}$-form) -- 2.4: \ } 
$(B_{\cal R}(u, s), \overline{r_u}, \overline{v}) 
\stackrel{2.4}{\longleftarrow} \ (s,  \overline{u}, \overline{v}) \ 
\stackrel{1.1}{\longrightarrow} \ (s, \overline{vu})$, \ \ \ 
where \ $s \leq_{\cal L} u$.

\medskip

{\bf Case 1.} \ $s \leq_{\cal L} vu \leq_{\cal L} u$.

Then \ $(s, \overline{vu}) \ \stackrel{2.4}{\longrightarrow} \ 
(B_{\cal R}(vu, s), \overline{r_{vu}})$.

On the other hand, \ $(B_{\cal R}(u, s), \overline{r_u}, \overline{v}) \  
\stackrel{1.1}{\longrightarrow} \ (B_{\cal R}(u, s), \overline{vr_u}) \ 
\stackrel{2.4}{\longrightarrow} \ (B_{\cal R}(vr_u, B_{\cal R}(u, s)),
\overline{r_{vr_u}})$. \ The last application of rule 2.4 is justified by 
Lemma \ref{Lm7} (1). 

To check confluence we observe that $vu \equiv_{\cal R} vr_u$, 
and that 
\ $B_{\cal R}(vu, s) = B_{\cal R}(vr_u, B_{\cal R}(u, s))$ \ by
Lemma \ref{Lm7} (1). 

\medskip

{\bf Case 2.} \ $vu <_{\cal L} s \leq_{\cal L} u$.

Then \ $(s, \overline{vu}) \ \stackrel{2.1}{\longrightarrow} \ 
(r_s, \overline{B_{\cal R}(s, vu)})$.

On the other hand, \ $(B_{\cal R}(u, s), \overline{r_u}, \overline{v}) \  
\stackrel{1.1}{\longrightarrow} \ (B_{\cal R}(u, s), \overline{vr_u}) \ 
\stackrel{2.1}{\longrightarrow} \ (r_{B_{\cal R}(u, s)}, 
\overline{B_{\cal R}(B_{\cal R}(u, s), vr_u)})$. \ 
The last application of rule 2.1 is justified by  Lemma  \ref{Lm7} (2).

Confluence now follows from Lemma  \ref{Lm7} (2), and from the fact that
$s \equiv_{\cal R} B_{\cal R}(u, s)$ (Lemma \ref{Lm2}).

\medskip

{\bf Case 3.} \ $vu \ \frac{<}{>}\!\!\!|\,_{_{\cal L}} \ s$.

Then \  $(s, \overline{vu}) \ \stackrel{1.3}{\longrightarrow} \ (0)$. \ 
On the other hand, \ $(B_{\cal R}(u, s), \overline{r_u}, \overline{v}) \
\stackrel{1.1}{\longrightarrow} \ (B_{\cal R}(u, s), \overline{vr_u}) \
\stackrel{1.3}{\longrightarrow} \ (0)$. \ We used Lemma \ref{Lm7} (3) to 
justify the last application of rule 1.3.

\bigskip

So far we have considered all overlaps involving the rule 1.1. We mentioned 
already that the rule 1.2 always leads to confluence to (0). Let us now look
at all the overlaps that involve rule 1.3 (other than with rule 1.1, 
seen already).

There is no overlap of 1.3 with itself.

\bigskip

\noindent {\bf Overlap 1.3--1.4: \ }  $(0, s) \stackrel{1.3}{\longleftarrow} 
\ (u, \overline{v}, s) \ \stackrel{1.4}{\longrightarrow} \  (u, 0)$, \ \ \    
where \ $u \ \frac{<}{>}\!\!\!|\,_{_{\cal L}} \ v$ \ and \ 
$v \ \frac{<}{>}\!\!\!|\,_{_{\cal R}} \ s$.\\ 
Then we obviously have confluence to (0).

The case of $(\overline{u}, v, \overline{s})$, \ where \ 
$u \ \frac{<}{>}\!\!\!|\,_{_{\cal R}} \ v$ \ and \  
$v \ \frac{<}{>}\!\!\!|\,_{_{\cal L}} \ v$, is handled in a similar way.

\bigskip

\noindent {\bf Overlap 1.3--1.5: \ }  $(B(u, v, w), \overline{s}) \  
\stackrel{1.5}{\longleftarrow} \ (u, \overline{v}, w, \overline{s}) \  
\stackrel{1.3}{\longrightarrow} \  (u, \overline{v}, 0)$, \\ 
where \  
$u \leq_{_{\cal L}} v \geq_{_{\cal R}} w$ \  and \  
$w \ \frac{<}{>}\!\!\!|\,_{_{\cal L}} \ s$.  

Then \ $(u, \overline{v}, 0) \longrightarrow (0)$ \ by two applications of 
rule 1.2. Moreover, since \ $B(u, v, w) \ \frac{<}{>}\!\!\!|\,_{_{\cal L}} \ s$
if \ $w \ \frac{<}{>}\!\!\!|\,_{_{\cal L}} \ s$ \ (by Lemma \ref{Lm8}), we
also have \ $(B(u, v, w), \overline{s}) \ 
\stackrel{1.3}{\longrightarrow} (0)$.  

\bigskip

\noindent {\bf Overlap 1.3--1.6: \ }  This is similar to 1.3--1.5.

\bigskip

\noindent There are no overlaps 1.3--2.1, 1.3--2.4, nor 1.4--1.4, 1.4--2.2, 
1.4--2.3. 
The overlaps 1.4--1.5  and 1.4--1.6 are similar to the case 1.3--1.5. 

\bigskip

\noindent {\bf Overlaps 1.3--2.2, 1.3--2.3, or 1.4--2.1: \ }  
This is very similar to the case considered next.

\bigskip

\noindent {\bf Overlap 1.4--2.4: \ } $(B_{_{\cal R}}(u, v), \overline{u}, w)
\ \stackrel{2.4}{\longleftarrow} \ (v, \overline{u}, w) \ 
\stackrel{1.4}{\longrightarrow} \ (v, 0)$, \ \ \  
where \ $v \leq_{_{\cal L}} u \ \frac{<}{>}\!\!\!|\,_{_{\cal R}} \ w$.

Then $(v, 0) \to (0)$ by rule 1.2. Moreover, since
$r_u \equiv_{_{\cal R}}  u \ \frac{<}{>}\!\!\!|\,_{_{\cal R}} \ w$ \ we 
have \ $(B_{_{\cal R}}(u, v), \overline{u}, w) \  \longrightarrow \ 
(B_{_{\cal R}}(u, v), 0)$ \ by rule 1.4; this then leads to (0) by 1.2. 

\bigskip

\noindent {\bf Overlap 1.5--1.5: \ } $(B(u, v, w), \overline{s}, t)
 \ \stackrel{1.5}{\longleftarrow} \  (u, \overline{v}, w, \overline{s}, t) \  
\stackrel{1.5}{\longrightarrow} \ (u, \overline{v}, B(w, s, t))$, \\
where \
$u \leq_{_{\cal L}} v \geq_{_{\cal R}} w \leq_{_{\cal L}} s \geq_{_{\cal R}} t$.

\smallskip

Then \ $(B(u, v, w), \overline{s}, t) \ \stackrel{1.5}{\longrightarrow} \ 
(B(B(u, v, w), s, t))$; \ rule 1.5 was applicable here by Lemma \ref{Lm9}. 
Also, \ $(u, \overline{v}, B(w, s, t)) \  
\stackrel{1.5}{\longrightarrow} \ (B(u, v, B(w, s, t))$; \ rule 1.5 was 
applicable here by Lemma \ref{Lm9}.
Confluence then follows from Lemma \ref{Lm9}.

\bigskip

\noindent {\bf Overlap 1.5--1.6: \ } $(B(u, v, w), \overline{s}) \  
\stackrel{1.5}{\longleftarrow} \ (u, \overline{v}, w, \overline{s}) \ 
\stackrel{1.6}{\longrightarrow} \ (u, \overline{B(s,w,v)})$, \  where
\ $u \leq_{_{\cal L}} v \equiv_{_{\cal R}} w \geq_{_{\cal L}} s$.

\medskip  

{\bf Case 1. \ } $B(u, v, w) \leq_{_{\cal L}} s$.

In this case rule 2.4 applies and 
\ $(B(u, v, w), \overline{s}) \ \stackrel{2.4}{\longrightarrow} \ 
(B_{_{\cal R}}(s, B(u, v, w)), \overline{r_s})$. \ 
By Lemma \ref{Lm11} (1), rule 2.4 then also applies to \ 
$(u, \overline{B(s, w, v)})$, \ thus producing \  
$(B_{_{\cal R}}(B(s, w, v), u), \overline{r_{B(s, w, v)}})$.
Lemma \ref{Lm11} (2.$\leq$) then shows confluence.

\medskip 

{\bf Case 2. \ } $B(u, v, w) >_{_{\cal L}} s$.

In this case $(B(u, v, w), \overline{s}) \ \stackrel{2.1}{\longrightarrow} \ 
(r_{B(u,v,w)}, \overline{B_{_{\cal R}}(B(u, v, w), s)})$. \ 
By Lemma \ref{Lm11} (1), rule 2.1 then also applies to \ 
$(u, \overline{B(s,w,v)})$, \ and this yields \ 
$(r_u, \overline{B_{_{\cal R}}(u, B(s,w,v))})$. 
Lemma \ref{Lm11} (2.$>$) then shows confluence.

\medskip 

{\bf Case 3. \ } $B(u, v, w) \ \frac{<}{>}\!\!\!|\,_{_{\cal L}} \ s$.

Then \ $(B(u, v, w), \overline{s}) \ \stackrel{1.3}{\longrightarrow} \ (0)$. 
Moreover, by Lemma \ref{Lm11} (1), in this case we also have \ 
$u \ \frac{<}{>}\!\!\!|\,_{_{\cal L}} \ B(s, w, v)$, \ hence rule 1.3 also 
applies to \ $(u, \overline{B(s,w,v)})$ \ and produces (0).

\medskip

\noindent The overlap case \ $\ \stackrel{1.6}{\longleftarrow} \ 
(\overline{u}, v, \overline{w}, s) \ \stackrel{1.5}{\longrightarrow} \ $
\ is similar to the case above.

\bigskip
\noindent {\bf Overlap 1.5--2.1: \ } $(B(u, v, w), \overline{s}) \ 
\stackrel{1.5}{\longleftarrow} \ (u, \overline{v}, w, \overline{s}) \
\stackrel{2.1}{\longrightarrow} \ 
(u, \overline{v}, r_w, \overline{B_{_{\cal R}}(w,s)})$,  \\  
where \  $u \leq_{_{\cal L}} v \geq_{_{\cal R}} w >_{_{\cal L}} s$. 

\medskip

{\bf Case 1. \ } $B(u, v, w) \leq_{_{\cal L}} s$.

Then \ $(B(u, v, w), \overline{s}) \ \stackrel{2.4}{\longrightarrow} \  
(B_{_{\cal R}}(s, B(u, v, w)), \overline{r_s})$. \ 
Moreover, \ $(u, \overline{v}, r_w, \overline{B_{_{\cal R}}(w,s)}) \  
\stackrel{1.5}{\longrightarrow} \ 
(B(u, v, r_w), \overline{B_{_{\cal R}}(w,s)}) \  
\stackrel{2.4}{\longrightarrow} \ 
(B_{_{\cal R}}(B_{_{\cal R}}(w,s), B(u, v, r_w)), 
\overline{r_{B_{_{\cal R}}(w,s)}})$. \ The last application of rule 2.4 is 
justified by Lemma \ref{Lm13} (1).
Confluence then follows immediately from Lemma \ref{Lm13} (2. $\leq$).

\medskip

{\bf Case 2. \ } $B(u, v, w) >_{_{\cal L}} s$.

Then \ $(B(u, v, w), \overline{s}) \ \stackrel{2.1}{\longrightarrow} \  
(r_{B(u, v, w)}, \overline{B_{_{\cal R}}(B(u, v, w), s)})$. \ 
Moreover, \ $(u, \overline{v}, r_w, \overline{B_{_{\cal R}}(w,s)}) \  
\stackrel{1.5}{\longrightarrow} \ 
(B(u, v, r_w), \overline{B_{_{\cal R}}(w,s)}) \  
\stackrel{2.1}{\longrightarrow} \ 
(r_{B(u,v,r_w)}, \overline{B_{_{\cal R}}(B(u,v,r_w), B_{_{\cal R}}(w,s))})$. 
\ The last application of rule 2.1 is justified by Lemma \ref{Lm13} (1).
Confluence then follows immediately from Lemma \ref{Lm13} (2. $<$).

\medskip

{\bf Case 3. \ } $B(u, v, w) \ \frac{<}{>}\!\!\!|\,_{_{\cal L}} \ s$.

Then \ $(B(u, v, w), \overline{s}) \ \stackrel{1.3}{\longrightarrow} \ (0)$.
Moreover, \ $(u, \overline{v}, r_w, \overline{B_{_{\cal R}}(w,s)}) \  
\stackrel{1.5}{\longrightarrow} \ 
(B(u, v, r_w), \overline{B_{_{\cal R}}(w,s)}) \  
\stackrel{1.3}{\longrightarrow} \  (0)$. \ The last application of rule 1.3
is justified by Lemma \ref{Lm13} (1).  
   
\bigskip

\noindent {\bf Overlap 1.5--2.2: \ } 

{\bf Case 1. \ } $u \leq_{_{\cal L}} v >_{_{\cal R}} w$ \ and \\  
$(B(u, v, w)) \ \stackrel{1.5}{\longleftarrow} \ 
(u, \overline{v}, w) \
\stackrel{2.2}{\longrightarrow} \ (u, \overline{\ell_v}, B_{_{\cal L}}(w,v)) \ 
\stackrel{1.5}{\longrightarrow} \ 
(B(u, \overline{\ell_v}, B_{_{\cal L}}(w,v)))$. 

Confluence then follows from the $B_{_{\cal L}}$-version of Lemma 
\ref{Lm14}.

\medskip

{\bf Case 2. \ }  $s >_{_{\cal R}} u \leq_{_{\cal L}} v \geq_{_{\cal R}} w$ \ 
and \\  
$(\overline{\ell_s}, B_{_{\cal L}}(u,s), \overline{v}, w) \  
\stackrel{2.2}{\longleftarrow} \ 
(\overline{s}, u, \overline{v}, w) \ \stackrel{1.5}{\longrightarrow} \ 
(\overline{s}, B(u, v, w))$. 

Then \ $(\overline{\ell_s}, B_{_{\cal L}}(u,s), \overline{v}, w) \
\stackrel{1.5}{\longrightarrow} \ 
(\overline{\ell_s}, B(B_{_{\cal L}}(u,s), v, w))$; \ rule 1.5 is applicable 
here since by Lemma \ref{Lm3}, \ $B_{_{\cal L}}(u,s) \equiv_{_{\cal L}} u 
\leq_{_{\cal L}} v \geq_{_{\cal R}} w$. 

On the other hand, \ 
$(\overline{s}, B(u, v, w)) \ \stackrel{2.2}{\longrightarrow} \ 
(\overline{\ell_s}, B_{_{\cal L}}(B(u, v, w), s))$; \ rule 2.2 is applicable 
here since \ $s >_{_{\cal R}} u \geq_{_{\cal R}} ux = B(u, v, w)$ \ 
(where the last equality holds by Lemma \ref{Lm5}).

Confluence then follows from the $B_{_{\cal L}}$-version of Lemma \ref{Lm15}. 

\bigskip

\noindent {\bf Overlap 1.5--2.3: \ } 

{\bf Case A.} \ \ $\ \stackrel{1.5}{\longleftarrow} \ 
(u, \overline{v}, w) \ \stackrel{2.3}{\longrightarrow} \ $, \ \ \ where \ 
$u \leq_{_{\cal L}} v \geq_{_{\cal R}} w$.

This is similar to Case A of the overlap 1.5--2.4, treated below.

\medskip

{\bf Case B.} \ \ $(\overline{B_{_{\cal L}}(s,u)}, \ell_u, \overline{v}, w) 
\ \stackrel{2.3}{\longleftarrow} \ 
(\overline{s}, u, \overline{v}, w) \   
\stackrel{1.5}{\longrightarrow} \ 
(\overline{s}, B(u, v, w))$, \ \ \ where \ 
$s \leq_{_{\cal R}} u \leq_{_{\cal L}} v \geq_{_{\cal R}} w$.
 
Then \ $(\overline{B_{_{\cal L}}(s,u)}, \ell_u, \overline{v}, w) \  
\stackrel{1.5}{\longrightarrow} \ 
(\overline{B_{_{\cal L}}(s,u)}, B(\ell_u, v, w))$.  

\smallskip

{\bf Case B.1} \ \ $s \leq_{_{\cal R}} B(u, v, w)$. 

Then \ $(\overline{s}, B(u, v, w)) \ \stackrel{2.3}{\longrightarrow} \  
(\overline{B_{_{\cal L}}(s, B(u, v, w))}, 
\ell_{B(u, v, w)})$.

On the other hand, \ $(\overline{B_{_{\cal L}}(s,u)}, B(\ell_u, v, w)) \ 
\stackrel{2.3}{\longrightarrow} \ 
(\overline{B_{_{\cal L}}(B_{_{\cal L}}(s,u), B(\ell_u, v, w))}, 
\ell_{B(\ell_u, v, w)})$.  \  Rule 2.3 was applicable here by the 
${\cal R}$-version of Lemma \ref{Lm13} (1).

Confluence then follows from the ${\cal R}$-version of 
Lemma \ref{Lm13} (2, $\leq$).

\smallskip

{\bf Case B.2} \ \ $s >_{_{\cal R}} B(u, v, w)$. 

Then \  $(\overline{s}, B(u, v, w)) \ \stackrel{2.2}{\longrightarrow} \  
(\overline{\ell_s}, B_{_{\cal L}}(B(u,v,w), s))$, \ and \\  
$(\overline{B_{_{\cal L}}(s,u)}, B(\ell_u, v, w)) \ 
\stackrel{2.2}{\longrightarrow} \ 
(\overline{\ell_{B_{_{\cal L}}(s,u)}}, 
B_{_{\cal L}}(B(\ell_u, v, w), B_{_{\cal L}}(s,u)))$. \ Rule 2.2 was 
applicable here by the ${\cal R}$-version of Lemma \ref{Lm13} (1). 

Confluence then follows from the ${\cal R}$-version of 
Lemma \ref{Lm13} (2, $<$).

\smallskip

{\bf Case B.3} \ \ $s \ \frac{<}{>}\!\!\!|\,_{_{\cal R}} \ B(u,v,w)$.

Then \  
$(\overline{s}, B(u, v, w)) \ \stackrel{1.4}{\longrightarrow} \ (0)$ 
\ and \\ 
$((\overline{B_{_{\cal L}}(s,u)}, B(\ell_u, v, w)) \
\stackrel{1.4}{\longrightarrow} \ (0)$, \ where the application of rule 1.4 is 
justified by the ${\cal R}$-version of Lemma \ref{Lm13} (1). 

\bigskip

\noindent {\bf Overlap 1.5--2.4: \ }  

\smallskip 

{\bf Case A.} \ $(B_{_{\cal R}}(v,u), \overline{r_v}, w) \ 
\stackrel{2.4}{\longleftarrow} \ (u, \overline{v}, w) \ 
\stackrel{1.5}{\longrightarrow} \ (B(u,v,w))$,  \\ 
where \ $u \leq_{_{\cal L}} v \geq_{_{\cal R}} w$.

Then rule 1.5 is applicable to \ $(B_{_{\cal R}}(v,u), \overline{r_v}, w)$ \  
because \ $u \leq_{_{\cal L}} v \geq_{_{\cal R}} w$ \ implies by Lemma 
\ref{Lm4} \ $B_{_{\cal R}}(v,u) \leq_{_{\cal L}} r_v \equiv_{_{\cal R}} v 
\geq_{_{\cal R}} w$. Applying 1.5 then yields 
$(B(B_{_{\cal R}}(v,u), r_v, w))$. \ Thus by Lemma \ref{Lm14} we have 
confluence.

\medskip

{\bf Case B.} \ $(B(u,v,w), \overline{s}) \ 
\stackrel{1.5}{\longleftarrow} \ (u, \overline{v}, w, \overline{s}) \ 
\stackrel{2.4}{\longrightarrow} \ 
(u, \overline{v}, B_{_{\cal R}}(s,w), \overline{r_s})$, \ \ 
where \ $u \leq_{_{\cal L}} v \geq_{_{\cal R}} w \leq_{_{\cal L}} s$.

Then rule 2.4 is applicable to \ $(B(u,v,w), \overline{s})$ \ because by Lemma
\ref{Lm5} \ $B(u,v,w) = yw \leq_{_{\cal L}} w \leq_{_{\cal L}} v$. Then 2.4
yields \ $(B_{_{\cal R}}(s, B(u,v,w)), \overline{r_s})$. 

On the other hand, rule 1.5 is applicable to \ 
$(u, \overline{v}, B_{_{\cal R}}(s,w), \overline{r_s})$ \ because \ 
$v \geq_{_{\cal R}} w \equiv_{_{\cal R}} B_{_{\cal R}}(s,w)$ \ 
(the latter by Lemma \ref{Lm3}). Then 1.5 yields \ 
$(B(u,v, B_{_{\cal R}}(s,w)), \overline{r_s})$. 

By Lemma \ref{Lm15} we have confluence. 

\bigskip

The overlaps of rule (1.6) with rules (1.6), (2.1)--(2.4) are handled in 
a similar way as the overlaps of (1.5) with rules (1.5), (2.1)--(2.4).

\bigskip

We now come to the overlaps of the rules 2.$i$ ($i$ = 1, ..., 4). 

Obviously, 2.1 cannot overlap with itself nor with 2.4.

\bigskip

\noindent {\bf Overlap 2.1--2.2: \ } $(r_u, \overline{B_{_{\cal R}}(u,v)}, w)
\ \stackrel{2.1}{\longleftarrow} \ (u, \overline{v}, w) \  
\stackrel{2.2}{\longrightarrow} \ 
(u, \overline{\ell_v}, B_{_{\cal L}}(w,v))$, \ \ 
where \ $u >_{_{\cal L}} v >_{_{\cal R}} w$. 

Then \ $(r_u, \overline{B_{_{\cal R}}(u,v)}, w) \ 
\stackrel{2.2}{\longrightarrow} \ 
(r_u, \overline{\ell_{B_{\cal R}(u,v)}}, B_{_{\cal L}}(w, B_{_{\cal R}}(u,v)))$.
\ Rule 2.2 was applicable here since by Lemma \ref{Lm3}, \ 
$B_{_{\cal R}}(u,v) \equiv_{_{\cal R}} v >_{_{\cal R}} w$. 

On the other hand, \  $(u, \overline{\ell_v}, B_{_{\cal L}}(w,v)) \  
\stackrel{2.1}{\longrightarrow} \ 
(r_u, \overline{B_{_{\cal R}}(u, \ell_v)}, B_{_{\cal L}}(w,v))$. \  
Rule 2.1 was applicable here since \  
$u >_{_{\cal L}} v \equiv_{_{\cal L}} \ell_v$. \ 

Next, applying rule 2.2 to this yields \  
$(r_u, \overline{\ell_{B_{\cal R}(u,\ell_v)}}, 
B_{_{\cal L}}(B_{_{\cal L}}(w,v), B_{_{\cal R}}(u,\ell_v)))$. 
\ Rule 2.2 was indeed applicable here since by Lemma \ref{Lm3}, 
\ $B_{_{\cal R}}(u,\ell_v) \equiv_{_{\cal R}} \ell_v \geq_{_{\cal R}} 
\ell_v y = B_{_{\cal L}}(w,v)$ \ where $uy = v$;
\ moreover, the $\geq_{_{\cal R}}$ is actually 
$>_{_{\cal R}}$ \ (if we had $\ell_v \equiv_{_{\cal R}} \ell_v y$, then we 
would also have $v \equiv_{_{\cal R}} vy = u$, which contradicts an 
assumption). 

Lemma \ref{Lm16} immediately shows confluence now.

\medskip

The other {\bf overlap case} for rules 2.1 and 2.2 is of the form

$(\overline{\ell_v}, B_{_{\cal L}}(v, w), \overline{w}) \  
\stackrel{2.2}{\longleftarrow} \ (\overline{u}, v, \overline{w}) \  
\stackrel{2.1}{\longrightarrow} \
(\overline{u}, r_v, \overline{B_{_{\cal R}}(v,w)})$, \\   
where \  $u >_{_{\cal R}} v >_{_{\cal L}} w$.

This case is similar to the case above.

\bigskip

\noindent {\bf Overlap 2.1--2.3: \ } 
$(r_u, \overline{B_{_{\cal R}}(u,v)}, w)
\ \stackrel{2.1}{\longleftarrow} \ (u, \overline{v}, w) \
\stackrel{2.3}{\longrightarrow} \
(u, \overline{B_{_{\cal L}}(v, w)}, \ell_w)$, \\   
where \  $u >_{_{\cal L}} v \leq_{_{\cal R}} w$.

Then \ $(r_u, \overline{B_{_{\cal R}}(u,v)}, w) \ 
\stackrel{2.3}{\longrightarrow} \
(r_u, \overline{B_{_{\cal L}}(B_{_{\cal R}}(u,v), w)}, \ell_w)$. \ 
Rule 2.3 was applicable here since \ 
$B_{_{\cal R}}(u,v) \equiv_{_{\cal R}} v$.

On the other hand, \  
$(u, \overline{B_{_{\cal L}}(v, w)}, \ell_w) \  
\stackrel{2.1}{\longrightarrow} \
(r_u, \overline{B_{_{\cal R}}(u, B_{_{\cal L}}(v,w))}, \ell_w)$. \
Rule 2.1 was applicable here since \
$B_{_{\cal L}}(v, w) \equiv_{_{\cal L}} v$. 

Confluence now follows from Lemma \ref{Lm17}.

\medskip

The other {\bf overlap case} for the rules 2.1 and 2.3 is of the form

$(\overline{B_{_{\cal L}}(u,v)}, \ell_v, \overline{w}) \   
\stackrel{2.3}{\longleftarrow} \ 
(\overline{u}, v, \overline{w}) \ 
\stackrel{2.1}{\longrightarrow} \
(\overline{u}, r_v, \overline{B_{_{\cal R}}(v,w)})$, \\   
where \ $u \leq_{_{\cal R}} v >_{_{\cal L}} w$.

This is similar to the overlap case of 2.2--2.4 that we will study next.

\bigskip

\noindent {\bf Rule 2.2} has no overlap with itself nor with 2.3.

\bigskip

\noindent {\bf Overlap 2.2--2.4: \ }  
$(B_{_{\cal R}}(v,u), \overline{r_v}, w) \  \stackrel{2.4}{\longleftarrow} \  
(u, \overline{v}, w) \  \stackrel{2.2}{\longrightarrow} \
(u, \overline{\ell_v}, B_{_{\cal L}}(w,v))$, \\
where \ $u \leq_{_{\cal L}} v >_{_{\cal R}} w$. 

Rule 1.5 is applicable to \ $(B_{_{\cal R}}(v,u), \overline{r_v}, w)$ \ since 
\ $B_{_{\cal R}}(v,u) = xr_v \leq_{_{\cal L}} r_v \equiv_{_{\cal R}} v 
>_{_{\cal R}} w$. \ This yields \ 
$(B(B_{_{\cal R}}(v,u), r_v, w))$. 

Rule 1.5 is also applicable to \ $(u, \overline{\ell_v}, B_{_{\cal L}}(w,v))$ \ 
since \ $u \leq_{_{\cal L}} v \equiv_{_{\cal L}} \ell_v \geq_{_{\cal R}} 
\ell_v y = B_{_{\cal L}}(w,v)$.  \ This yields \ 
$(B(u, \ell_v, B_{_{\cal L}}(w,v))$. 

Lemma \ref{Lm18} immediately implies confluence.

\medskip

The other {\bf overlap case} for the rules 2.2 and 2.4 is of the form

$(\overline{\ell_u}, B_{_{\cal L}}(v,u), \overline{w}) \
\stackrel{2.2}{\longleftarrow} \ (\overline{u}, v, \overline{w}) \
\stackrel{2.4}{\longrightarrow} \
(\overline{u}, B_{_{\cal R}}(w,v), \overline{r_w})$, \ \ \ where \  
$u >_{_{\cal R}} v \leq_{_{\cal L}} w$. 

This is very similar to the overlap case of 2.1--2.3 that we studied 
explicitly.

\bigskip

\noindent {\bf Overlap 2.3--2.4: \ }
$(B_{_{\cal R}}(v,u), \overline{r_v}, w) \  \stackrel{2.4}{\longleftarrow} \
(u, \overline{v}, w) \  \stackrel{2.3}{\longrightarrow} \
(u, \overline{B_{_{\cal L}}(v,w)}, \ell_w)$, \\  
where \ $u \leq_{_{\cal L}} v \geq_{_{\cal R}} w$.

Then \ $(B_{_{\cal R}}(v,u), \overline{r_v}, w) \  
\stackrel{2.3}{\longrightarrow} \
(B_{_{\cal R}}(v,u), \overline{B_{_{\cal L}}(r_v, w)}, \ell_w)
\stackrel{2.4}{\longrightarrow} \
(B_{_{\cal R}}(B_{_{\cal L}}(r_v, w), B_{_{\cal R}}(v,u)), 
 \overline{r_{B_{_{\cal L}}(r_v, w)}}, \ell_w)$; \ the last application of 
rule 2.4 was justified since \ $B_{_{\cal R}}(v,u) = xr_v \leq_{_{\cal L}}
r_v \equiv_{_{\cal L}} B_{_{\cal L}}(r_v, w)$ \ (the last ${\cal L}$-equivalence
follows from Lemma \ref{Lm3}). 

On the other hand, \ $(u, \overline{B_{_{\cal L}}(v,w)}, \ell_w) \ 
\stackrel{2.4}{\longrightarrow} \
(B_{_{\cal R}}(B_{_{\cal L}}(v,w), u), 
\overline{r_{B_{_{\cal L}}(v,w)}}, \ell_w)$; \ 
the application of rule 2.4 was
justified since \ $u \leq_{_{\cal L}} v \equiv_{_{\cal L}} B_{_{\cal L}}(v,w)$
\ (where the last ${\cal L}$-equivalence follows from Lemma \ref{Lm3}). 

Confluence now follows immediately from the ${\cal L}-{\cal R}$ dual of 
Lemma \ref{Lm16}.

\medskip 

The other {\bf overlap case} for the rules 2.3 and 2.4 is of the form

$(\overline{B_{_{\cal L}}(u,v)}, \ell_v, \overline{w}) \
\stackrel{2.3}{\longleftarrow} \
(\overline{u}, v, \overline{w}) \
\stackrel{2.4}{\longrightarrow} \
(\overline{u}, B_{_{\cal R}}(w,v), \overline{r_v})$, \\
where \ $u \leq_{_{\cal R}} v \geq_{_{\cal L}} w$.

This is similar to the above case.

\bigskip

This completes the exhaustive analysis of all overlap cases, and shows that the 
rewrite system for $(S)_{\mathrm{reg}}$ is {\it locally confluent}.


\begin{thebibliography}{99}

\bibitem{Bi1} J.C.\ Birget, ``Iteration of expansions -- unambiguous
semigroups'', {\it J.\ of Pure and Applied Algebra} 34 (1984) 1-55.

\bibitem{Bi2} J.C.\ Birget, ``Arbitrary vs.\ regular semigroups'',
{\it J.\ of Pure and Applied Algebra} 34 (1984) 57-115.

\bibitem{Bi3} J.C.\ Birget, ``Time-complexity of the word problem for
semigroups and the Higman embedding theorem'', 
{\it International J.\ of Algebra and Computation} 8 (1998) 235-294.

\bibitem{Bi4} J.C.\ Birget, ``Historical and Technical Perspective on the Synthesis
Theorem", in {\it Monoids and Semigroups with Applications} 
(J.\ Rhodes, editor), Proc.\ of  1989 Berkeley Workshop, World Scientific 
Publ.\ Co.\ (1991), pp.\ 393-402.

\bibitem{GrBook} P.\ A.\ Grillet, {\it Semigroups: An Introduction to the 
Structure Theory}, Marcel Dekker (1995).

\bibitem{GrPaper} P.\ A.\ Grillet, ``On Birget's regular embedding'', 
{\it J.\ of Pure and Applied Algebra} 130 (1998) 293-311. 

\bibitem{Ja} M.\ Jantzen, {\it Confluent String Rewriting}, Springer-Verlag
(1988).

\bibitem{Pecuchet} J.P.\ P\'ecuchet, ``Automates boustrophedon, semigroupe de
Birget et mono\"{\i}de inversif libre'', {\it RAIRO Informatique Th\'eorique}
19 (1985) 71-100.

\end{thebibliography}
\end{document}